\documentclass[12pt]{extarticle}

\usepackage{booktabs} % For better-looking tables

\usepackage{amsmath,amssymb,amsthm}
\usepackage{bbm}
\usepackage[margin=1in]{geometry}
\usepackage{graphicx}
\usepackage{hyperref}
\usepackage[utf8]{inputenc}
\usepackage[T1]{fontenc}
\usepackage{mathtools}
\usepackage[english]{babel}

\usepackage{placeins}
\usepackage{hyperref}
\hypersetup{
    colorlinks,
    linkcolor={blue!80!black},
    citecolor={blue},
    urlcolor={black}
}
\usepackage{xcolor}

\usepackage{xcolor}

\usepackage{enumitem}

\usepackage{float}

\usepackage{tikz}
\usetikzlibrary{arrows.meta, patterns}

\usepackage{cancel}

\usepackage{breqn}

\let\div\relax
\DeclareMathOperator{\div}{div}

\DeclareMathOperator{\Range}{Range}

\usepackage{lmodern}
\usepackage{fix-cm}

\DeclareMathSizes{14}{14}{9}{7}

\numberwithin{equation}{section}

\newcommand{\keywords}[1]{\par\addvspace{2\baselineskip}\noindent\textbf{Keywords:} #1}

\newcommand{\bH}{\mathbf{H}}
\newcommand{\bL}{\mathbf{L}}
\newcommand{\bu}{\mathbf{u}}
\newcommand{\bw}{\mathbf{w}}
\newcommand{\bv}{\mathbf{v}}
\newcommand{\bz}{\mathbf{z}}
\newcommand{\bX}{\mathbf{X}}
\newcommand{\bZ}{\mathbf{Z}}
\newcommand{\bY}{\mathbf{Y}}
\newcommand{\bV}{\mathbf{V}}
\newcommand{\bW}{\mathbf{W}}
\newcommand{\bG}{\mathbf{G}}
\newcommand{\bC}{\mathbf{C}}

\newcommand{\bg}{\mathbf{g}}
\newcommand{\bSigma}{\mathbf{\Sigma}}

\newcommand{\bpsi}{\boldsymbol{\psi}}
\newcommand{\bbeta}{\boldsymbol{\eta}}
\newcommand{\bomega}{\boldsymbol{\omega}}
\newcommand{\bvarphi}{\boldsymbol{\varphi}}

\newcommand{\btau}{\boldsymbol{\tau}}

\newcounter{statementcounter}

\newtheoremstyle{customnormal}% name
  {\topsep}
  {\topsep}
  {\normalfont}
  {}%         Indent amount (empty = no indent)
  {\bfseries}
  {.}%         Punctuation after theorem head (a period)
  {.5em}%         Space after theorem head
  {\thmname{#1}\thmnumber{ \thestatementcounter}\thmnote{ (\normalfont#3)}}

\theoremstyle{customnormal}
\newtheorem{mytheorem}[statementcounter]{Theorem}
\newtheorem{myproposition}[statementcounter]{Proposition}
\newtheorem{mylemma}[statementcounter]{Lemma}
\newtheorem{mycorollary}[statementcounter]{Corollary}

\title{A Revisiting of the Pressure Elimination for a Fluid--Structure PDE Interaction and Its Implications}

\author{
George Avalos\thanks{Department of Mathematics, University of Nebraska-Lincoln, Lincoln, NE, \href{mailto:gavalos2@unl.edu}{\texttt{gavalos2@unl.edu}}  }
\and 
Yuhao Mu\thanks{Department of Mathematics, University of Nebraska-Lincoln, Lincoln, NE, \href{mailto:ymu4@huskers.unl.edu}{\texttt{ymu4@huskers.unl.edu}} }
}
\date{}

\begin{document}
\maketitle

\begin{abstract}
In this paper we construct a novel technique for eliminating and recovering the pressure for a fluid--structure interaction model. This pressure elimination methodology is valid for general bounded Lipschitz domains. The specific fluid--structure interaction (FSI) that we consider is a well-known model of Stokes flow coupled to a system of linear elasticity, which constitutes a coupled parabolic-hyperbolic system. The coupling between the two distinct PDE dynamics occurs across a boundary interface, with each of the components evolving on its own distinct geometry, with the domains of each being Lipschitz. Our new pressure elimination technique admits of an explicit $C_{0}$-semigroup generator representation $\mathcal{A}: D(\mathcal{A}) \subset \mathbf{H} \to \mathbf{H}$, where $\mathbf{H}$ is the associated finite energy space of fluid--structure initial data. This leads to a novel proof of well-posedness in the explicit semigroup sense of the continuous PDE, now valid in general geometries. Subsequently, we illustrate an immediate consequence of our semigroup well-posedness result; namely a finite element method (FEM) with associated rates of convergence for a static version of the FSI, posed on polygonal domains.
\end{abstract}


\section{Introduction}
\begin{figure}[htbp]
  \centering
  \includegraphics[scale=1]{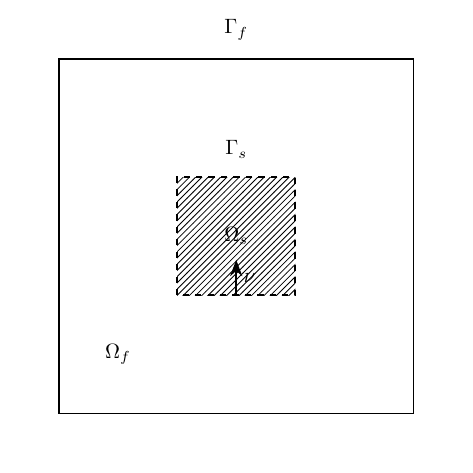}
  \vspace{-0.5cm} % Reduce space below
  \caption{The fluid--structure geometry.}
  \label{fig:fsi_geometry}
\end{figure}

We consider a model of fluid--structure interaction on two bounded, Lipschitz domains $\Omega_f$ and $\Omega_s$, where each $\Omega_i \subset \mathbb{R}^n$, $n = 2$ or $3$. (A canonical example of the pair $\{ \Omega_f, \Omega_s \}$ is given in Figure~\ref{fig:fsi_geometry}. The physical situation to be modeled is as follows: a stationary elastic solid $\Omega_s$ is fully immersed in a fluid occupying domain $\Omega_f$ with the interaction taking place on the boundary of the solid $\Gamma_s$. The dynamics of the solid is described by a linear elastic (hyperbolic) equation in the variable $\bw$, while the velocity of the fluid $\bu$ is modeled by (parabolic-like) Stokes flow.
\par To be specific, we consider the following linear version of a nonlinear fluid--structure interaction (FSI) model combining a solid displacement (linear elasticity) with Stokes flow:
\begingroup
  \setlength{\jot}{2pt}
  \begin{subequations}
  \begin{alignat}{2}
    \bu_t - \div(\varepsilon(\bu)) + \nabla p = 0
    &&\quad \text{in }\Omega_f\times(0,T),
      \label{eq:fsif1}\\
    \div(\bu) = 0
    &&\quad \text{in }\Omega_f\times(0,T),
      \label{eq:fsif2}\\
    \bw_{tt} - \div(\sigma(\bw)) + \bw = 0
    &&\quad \text{in }\Omega_s\times(0,T),
      \label{eq:fsif3}\\
    \bu = 0
    &&\quad \text{in }\Gamma_f\times(0,T),
      \label{eq:fsif4}\\
    \bw_t = \bu
    &&\quad \text{on }\Gamma_s\times(0,T),
      \label{eq:fsif5}\\
    \sigma(\bw)\cdot\nu = \varepsilon(\bu)\cdot\nu - p\,\nu
    &&\quad \text{on }\Gamma_s\times(0,T),
      \label{eq:fsif6}\\
    \bu(\cdot,0) = \bu_0
    &&\quad \text{in }\Omega_f,
      \label{eq:fsif7}\\
    \bw(\cdot,0)=\bw_0,\;\bw_t(\cdot,0)=\bw_1
    &&\quad \text{in }\Omega_s.
      \label{eq:fsif8}
  \end{alignat}
  \end{subequations}
\endgroup
So, as presented, we have $\Omega_f$ is the fluid region and $\Omega_s$ the structure region of the FSI. Also, the stress tensor $\sigma(\cdot)$ and strain tensor $\varepsilon(\cdot)$ are defined as in \eqref{eq:stress_tensor} and \eqref{eq:strain_tensor}. The $\Gamma_f$ region is the outer boundary of the fluid (for example, the bigger circle in an annulus), and $\Gamma_s$ is the inner boundary of the structure (or between the structure and the fluid). The interchange between the fluid and the structure occurs at this boundary, as written above. By our convention $\nu(\mathbf{x})$ is the unit outward normal vector with respect to $\Omega_f$ and hence inward with respect to $\Omega_s$.

\par This is the linear version of a nonlinear system that is well-known in the literature (see \cite{barbu2007existence}\cite{barbu2008smoothness}) which uses Navier--Stokes rather than Stokes flow, with the nonlinear (advection) term and nonlinear coupling term removed. For other treatments of this system see, for e.g.: \cite{avalos2007semigroup}, \cite{avalos2009semigroup}, \cite{avalos2008maximality}, \cite{raymond2014fsi}. For the nonlinear analysis of this type of system, there is also, for e.g.: \cite{desjardins2001weak}, \cite{grandmont2002existence}, \cite{beirao2004strong}, \cite{bociu2011linearization}, \cite{chueshov2013unsteady}, and see also the references in Section~\ref{past_techniques}.

\vspace{0.2cm}

\textbf{Outline and goals:} The main goal of this paper is to introduce a new method of eliminating the pressure $p(t)$ from the FSI system in Section~\ref{newmethod}, and in conjunction with the Lumer--Phillips theorem (see \cite{pazy1983semigroups}), apply it to obtain semigroup well-posedness of the FSI modeled \eqref{eq:fsif1}--\eqref{eq:fsif8} (see Section~\ref{FEM_estimates}).
Subsequently, we provide a companion finite element method (FEM), with explicit rates of convergence, to approximate solutions of the associated FSI semigroup resolvent equations.

We emphasize that the use of our pressure elimination method is not limited to the linear fluid-elasticity model \eqref{eq:fsif1}--\eqref{eq:fsif8}, but in principle can be extended to a range of other (both linear and nonlinear) fluid--structure and poroelastic-fluid flow models, under the various modes of interactive boundary coupling (see \cite{avalos2025poroelastic}).

Our present approach to eliminate the associated pressure for the FSI model \eqref{eq:fsif1}--\eqref{eq:fsif8} is altogether different than that invoked in the analysis of other fluid--structure and Biot--Stokes systems; see, e.g. \cite{avalos2007semigroup}, \cite{avalos2009semigroup}, \cite{avalos2008maximality}, \cite{avalos2014mixed}, \cite{avalos2015rational}, \cite{geredeli2025_spectral}, \cite{bociu2015linearization}, \cite{avalos2025poroelastic}, \cite{sviridyuk2022avalostr}. Namely, in these works, the associated fluid pressure variable is characterized as the solution of a certain elliptic boundary value problem (BVP) (see \ref{past_techniques}).

This (initially formal) BVP-characterization of the pressure allows for a proper description of the matrix operator $\mathcal{A}$, say, this (candidate semigroup generator) mathematically models the various differential actions and boundary coupling mechanisms for the given FSI model.

However, in order to properly define the domain of the modeling operator $\mathcal{A}$, in such a way as to ultimately give rise to its strongly continuous semigroup generation, the said BVP-characterization of the pressure variable has to be made wholly rigorous.

In particular, some measure of elliptic regularity must already be in hand, so as to insure that (negative) Sobolev regularity of certain fluid and structure boundary traces, as they appear as a boundary data of the pressure BVP---see, e.g., \ref{past_techniques} below---will give rise to the pressure variable, as the solution of the given boundary value problem (see the more detailed discussion in Section~\ref{past_techniques} below). But to have an available elliptic regularity theory will require smoothness assumptions on the fluid--structure geometries. Such smoothness assumption are made in the aforementioned works.

However, the geometric smoothness assumptions, which are needed to assure an underlying elliptic regularity, will preclude the consideration of even the canonical fluid--structure domain in Figure~\ref{fig:fsi_geometry} above. In particular, the non-convexity and obtuseness of the angles of fluid domain $\Omega_f$ (about the boundary interface $\Gamma_s$) will not allow for full $H^2$-regularity of second order elliptic BVP solutions, with homogeneous boundary conditions (see \cite{grisvard1985_elliptic}, \cite{dauge1988_corner}). This higher regularity is required in the ``Method of Transposition'', of J.-L. Lions and E. Magenes, by way of solving elliptic BVPs with weak boundary data (see \cite{lions1972nhbvp}).

Instead, in the present work, we eliminate the associated pressure variable by other means. As was the case for smooth domains, this novel pressure elimination will likewise allow for the construction of a modeling operator $\mathcal{A} : \bH \to \bH$ of \eqref{eq:fsif1}--\eqref{eq:fsif8}, for which there is $C_0$-semigroup generation. (There is still a BVP-characterization of pressure variable $p(t)$ of \eqref{eq:fsif1}--\eqref{eq:fsif8}, but for Lipschitz pair $\{\Omega_f, \Omega_s \}$, it is purely formal.)

Our new way of proceeding is explicitly formulated in Section~\ref{newmethod} below, but we can state the heart of the matter here: We will use the ``matching fluxes'' boundary condition in \eqref{eq:fsif6} to ultimately put ourselves in position to invoke the classic result of de Rham to characterize the ``zero mean average''-component of pressure variable $p(t)$ (see \cite{temam1977_navier} and \cite{sohr2012navier}). Putting ourselves in position ostensibly entails an appropriate and different description of the domain of the semigroup generator $D(\mathcal{A})$, vis-a-vis that constructed for sufficiently smooth geometries. However, $D(\mathcal{A})$ in the present nonsmooth case is \textit{exactly the same} as that for smooth geometries, save that there is now no need to impose (negative) Sobolev regularity of boundary trace $\Delta \bu \cdot \nu|_{\Gamma_f}$.

Having in hand FSI semigroup generation for Lipschitz pair $\{ \Omega_f, \Omega_s \}$, one can subsequently speak of FEM approximation of (initially) static versions of the system \eqref{eq:fsif1}--\eqref{eq:fsif8}. This is undertaken in Sections~\ref{past_techniques} and \ref{numerical_results}.

\vspace{0.2cm}

\textbf{Notation:} in the following, vector-valued functions and spaces of such functions (excepting the unit normal $\nu$) will be denoted with boldface. For example, we denote $\bH^{1}(\Omega) \equiv [H^{1}(\Omega)]^n$. We define the energy space:
\begin{equation*}\bH = \bH_{f} \times \bH^{1}(\Omega_s) \times \bL^{2}(\Omega_s),\end{equation*}
with the fluid component space being defined as:
\begin{equation*}\bH_{f}\equiv\{\bu \in \bL^{2}(\Omega_f) : \div (\bu) = 0, \bu \cdot \nu|_{\Gamma_f} \,=\, 0 \},\end{equation*}
where we are using the fact that if $\bu \in \bL^{2}(\Omega_f)$ and $\div(\bu) \in L^{2}(\Omega_f)$, one has $\bu \cdot \nu|_{\partial\Omega_f} \in H^{-1/2}(\partial\Omega_f)$, and so $\bH_f$ is well-defined (see \cite{constantin1988navier}). Recall also that the stress and strain tensors are:
\begin{align}
\sigma(\mathbf{w}) &= \lambda_{s}\,\mathrm{tr}\big(\varepsilon(\mathbf{w})\big)\,I
                     + 2\mu\,\varepsilon(\mathbf{w}),
    \label{eq:stress_tensor} \\[4pt]
\varepsilon(\mathbf{w}) &= \tfrac{1}{2}\big(\nabla \mathbf{w}
                        + (\nabla \mathbf{w})^{T}\big).
    \label{eq:strain_tensor}
\end{align}
We also have $\bH$ is a Hilbert space with the following norm inducing inner product:

\begin{equation}\left( \begin{bmatrix}\bu\\ \bw\\ \bz\end{bmatrix}, \begin{bmatrix}\tilde{\bu}\\ \tilde{\bw}\\ \tilde{\bz}\end{bmatrix} \right)_{H} = (\bu, \tilde{\bu})_{\Omega_f} + (\bw, \tilde{\bw})_{1, \Omega_s} + (\bz, \tilde{\bz})_{\Omega_s},\end{equation}
where $(\mathbf{f}, \mathbf{g})_{\Omega} = \int_{\Omega}\mathbf{f}\cdot\mathbf{g}\,d\Omega$, and:
\begin{equation}
(\bu, \tilde{\bu})_{1, \Omega_f} = (\varepsilon(\bu), \varepsilon(\tilde{\bu}))_{\Omega_f} \quad \text{ : } \quad (\bw, \tilde{\bw})_{1, \Omega_s} = (\varepsilon(\bw), \sigma(\tilde{\bw}))_{\Omega_s} + (\bw, \tilde{\bw})_{\Omega_s}.
\end{equation}
In addition, we define the space of fluid test functions:
\begin{equation}\label{test_fcns}
\bH_{\Gamma_f, 0}^{1}(\Omega_f) \equiv \{\bvarphi \in \bH^{1}(\Omega_f) : \bvarphi|_{\Gamma_f} = 0 \},\end{equation}
which is topologized with respect to the norm $|\cdot|_{1, \Omega_f}$ induced by $(\bvarphi, \tilde{\bvarphi})_{1, \Omega_f}$ which is equivalent to the usual $\bH^{1}(\Omega_f)$ norm via Korn's inequality and Poincar\'e's inequality. Similarly we have $\bH^{1}(\Omega_s)$ is topologized with respect to the norm $|\cdot|_{1, \Omega_s}$ given by $(\bw, \tilde{\bw})_{1, s}$ which is equivalent to the usual Sobolev $\bH^{1}(\Omega_s)$ norm by Korn's inequality (see \cite{kesavan1989topics}):
\begin{equation}
(\varepsilon(\bw), \varepsilon(\bw))_{\Omega_s} \,+\, (\bw, \bw)_{\Omega_s} \;\geq\; c \|\bw\|_{\bH^{1}(\Omega_s)}^2 \text{ where } c \geq 0,
\end{equation}
where in the following we use $\|\cdot\|_{k, \Omega}$ to denote the Sobolev norm $\|\cdot\|_{\bH^{k}(\Omega)}$ of order $k$ over domain $\Omega$.

\section{Elimination of the pressure}

\subsection{Previous techniques to eliminate the pressure}
\label{past_techniques}
\par This fluid--structure PDE appeared in \cite{lions1969quelques} and subsequently in \cite{du2003analysis}; there is also the work of \cite{avalos2007semigroup}, \cite{barbu2007existence},  \cite{barbu2008smoothness}, \cite{avalos2008maximality}, \cite{avalos2009semigroup}, \cite{avalos2014mixed}, \cite{ignatova2014nonlinearity}, \cite{avalos2015rational}, \cite{geredeli2025_spectral}, \cite{bociu2015linearization}, \cite{avalos2025poroelastic}, \cite{sviridyuk2022avalostr}. In these works, there is a generally recognized need to deal with the pressure (necessarily in a nonstandard fashion). Specifically in the the coupled problem \eqref{eq:fsif1}--\eqref{eq:fsif8}, due to the non-homogeneous boundary coupling, it is \textit{not} possible to use the classical approach for no-slip boundary conditions to eliminate the pressure: that is, by applying the Leray projector on the equation from $\bL^{2}(\Omega_f)$ onto the classical space, which supposes $\bu \in \bL^{2}(\Omega_f)$, $\div(\bu) = 0$ on $\Omega$, $\bu \cdot \nu = 0$ on $\partial\Omega_f$. 

\par The following is the approach which was taken in \cite{avalos2007semigroup} and \cite{avalos2008maximality}. For the sake of simplicity, we consider the linear (Stokes) version of the FSI. The associated pressure $p(t)$ necessarily satisfies the following elliptic problem:
\begin{equation}
\begin{aligned}
\Delta p &= 0 \quad \text{in } \Omega_f, \\
p &= (\varepsilon(\bu) \cdot \nu)\cdot \nu - (\sigma(\bw) \cdot \nu) \cdot \nu \quad \text{on } \Gamma_s, \\
\frac{\partial p}{\partial\nu} &= \text{div}(\varepsilon(\bu) \cdot \nu) \quad \text{on } \Gamma_f.
\end{aligned}
\end{equation}
Such characterization of the pressure for uncoupled flows has been long recognized, see e.g. \cite{doering1995navier}, where a characterizing BVP for the pressure is formally undertaken with respect to uncoupled Euler's equations.
\par As such, the pressure then admits of the formal representation:
\begin{equation}
p(t) = D_{s}\left[ (\varepsilon(\bu) \cdot \nu) \cdot \nu - (\sigma(\bw) \cdot \nu) \cdot \nu \right]_{\Gamma_s} 
+ N_{f} \left[\text{div}(\varepsilon(\bu)) \cdot \nu |_{\Gamma_f} \right] \quad \text{in } \Omega_f;
\end{equation}
where the ``Dirichlet'' map $D_s$ is defined by:
\begin{equation}
h = D_{s}(g) \iff 
\begin{cases} 
\Delta h = 0  & \text{in } \Omega_f, \\ 
h = g  & \text{on } \Gamma_s, \\ 
\frac{\partial h}{\partial \nu} = 0 & \text{on } \Gamma_f, 
\end{cases}
\end{equation}
and the ``Neumann'' map $N_f$ is the map defined by:
\begin{equation}
h = N_{f}(g) \iff 
\begin{cases} 
\Delta h = 0  & \text{in } \Omega_f, \\ 
h = 0  & \text{on } \Gamma_s, \\ 
\frac{\partial h}{\partial \nu} = g & \text{on } \Gamma_f. 
\end{cases}
\end{equation}
Upon substitution, the PDE fluid component \eqref{eq:fsif1} thus becomes
\begin{equation}
\bu_{t} = \text{div}\varepsilon(\bu) + \widetilde{\bG}_{1}\bw + \widetilde{\bG}_{2}\bw \quad \text{in } (0, T) \times \Omega_f,
\end{equation}
where we have the maps:
\begin{equation}
\begin{aligned}
\widetilde{\bG}_{1}\bw &\equiv \nabla \left(D_{s}\left[(\sigma(\bw) \cdot \nu) \cdot \nu|_{\Gamma_s} \right]\right), \\
\widetilde{\bG}_{2}\bu &\equiv -\nabla D_{s}\left[(( \varepsilon(\bu) \cdot \nu) \cdot \nu)|_{\Gamma_{s}} \right] 
-\nabla N_{f} \left[\div(\varepsilon(\bu)) \cdot \nu |_{\Gamma_f}\right].
\end{aligned}
\end{equation} 

This elimination of the associated FSI pressure via the BVP above initially proceeded by formal calculation. In the course of the previously cited work, this characterization of the pressure is made rigorous: Namely, given the $C^2$ smoothness of the pair $\{\Omega_s, \Omega_f \}$ one will have full $H^2(\Omega_f)$-regularity of Poisson's equation with square integrable forcing term and homogeneous BCs (see \cite{lions1972nhbvp}), viz., if $q$ solves 
\begin{equation}
-\Delta q = f \in L^{2}(\Omega_f), \quad q|_{\Gamma_s} = 0 \text{ on } \Gamma_s, \quad \frac{\partial q}{\partial\nu}\bigg|_{\Gamma_s} = 0, \quad q \in H^{2}(\Omega_f).
\end{equation}  
In turn, via the Transposition Method of J.-L. Lions and E. Magenes \cite{lions1972nhbvp} for finding solutions of the BVP above in negative Sobolev spaces, one can give a $L^2$-meaning to the solution $p(t)$ for data
\begin{equation}
\label{boundary_traces}
\begin{aligned}
&\left[(\varepsilon(\bu) \cdot \nu) \cdot \nu - (\sigma(\bw) \cdot \nu) \cdot \nu \right]_{\Gamma_s} \in H^{-1/2}(\Gamma_s); \\
&\left(\text{div}(\varepsilon(\bu)) \cdot \nu |_{\Gamma_f} \right) \in H^{-3/2}(\Gamma_f).
\end{aligned}
\end{equation}   
(Note that the Method of Transposition implicitly makes the further requirement that $\Gamma_f \in C^3$. See Theorem 2.7.4, p. 101, of \cite{kesavan1989topics} and Proposition 1, p. 265 of \cite{avalos2008maximality}.) Specifically the harmonic extension maps $D_{s}$ and $N_{f}$ into $L^{2}(\Omega_f)$ are well-defined if the corresponding Poisson problems with homogeneous Neumann or Dirichlet conditions define topological isomorphisms of $L^{2}(\Omega_f)$ onto the corresponding subspace of $H^{2}(\Omega_f)$. Such elliptic smoothing however, is known to break down even in the 2D polygonal setting (see for example Corollary 4.4.3.8 in \cite{grisvard1985_elliptic}, or \cite{dauge1988_corner}). In turn, the modeling $C_{0}$-semigroup generator $\mathcal{A}$, whose form is dictated by the substitution above, then has its domain specified so these ``smooth'' solutions — i.e., solutions corresponding to data in $D(\mathcal{A})$ — has the well-defined boundary traces in \eqref{boundary_traces}. This idea to eliminate the (incompressible) pressures associated with FSI and Biot--Stokes has been invoked in \cite{bociu2015linearization}, \cite{avalos2015rational}, \cite{avalos2014mixed}, \cite{sviridyuk2022avalostr}, \cite{avalos2025poroelastic}. However to make the said pressure substitution valid, the geometry has to be sufficiently smooth. This obstacle (and our following solution to it) has a number of consequences for computation: For example, \cite{avalos2016infsup} demonstrated that the mixed FEM above, with the so-called Taylor--Hood formulation in place for the fluid PDE component, obeys the discrete Babu\v{s}ka--Brezzi condition, uniformly, with respect to the discretization parameter. Thus, since our following solution makes the FSI generator rigorous for polygonal domains (and generally, domains for which the FEM can be applied), one now has FEM convergence estimates immediately at hand (see Section~\ref{FEM_estimates}).

\subsection{New method of pressure elimination}
\label{newmethod}
By way of obtaining a valid fluid--structure semigroup generator for general Lipschitz domain pair $\{ \Omega_f, \Omega_s\}$, we undertake in this section an elimination of the associated pressure variable in a manner which does not require any underlying elliptic regularity, in contrast to what was done in \cite{avalos2008maximality}, \cite{avalos2007semigroup} \cite{avalos2009semigroup} (those operating under sufficiently smooth geometric assumptions).
\par To this end, if $p(\mathbf{x}, t)$ in \eqref{eq:fsif1} is in $L^{2}(\Omega_f)$, a.e. in time: we first consider the fluid PDE component \eqref{eq:fsif1}. Multiplying both sides by $\bvarphi \in \bV$, where \begin{equation}\label{test_space}\bV\equiv \{\bvarphi \in \bH^{1}(\Omega_f) : \bvarphi|_{\Gamma_f} \,= 0,\; \div(\bvarphi) = 0 \text{ in } \Omega_f \},\end{equation} and integrating by parts, we get:
\begin{equation*}
\begin{split}
(\bu_t, \bvarphi)_{\Omega_f} &=
(\div(\varepsilon(\bu)), \bvarphi)_{\Omega_f} - (\nabla p, \bvarphi)_{\Omega_f}\\
& = \langle \varepsilon(\bu) \cdot \nu, \bvarphi\rangle_{\Gamma_s} - (\varepsilon(\bu), \varepsilon(\bvarphi))_{\Omega_f}
- \langle p, \bvarphi \cdot \nu\rangle_{\Gamma_s} + \cancel{(p, \div(\bvarphi))_{\Omega_f}} \\
& = \langle \varepsilon(\bu) \cdot\nu - p\nu, \bvarphi\rangle_{\Gamma_s} - (\varepsilon(\bu), \varepsilon(\bvarphi))_{\Omega_f}\\
& = \langle \sigma(\bw) \cdot\nu, \bvarphi\rangle_{\Gamma_s} - (\varepsilon(\bu), \varepsilon(\bvarphi))_{\Omega_f},
\end{split}
\end{equation*}
where in the last relation, we invoked the matching fluxes BC in \eqref{eq:fsif6}. Then by the Riesz representation theorem we have:
\begin{equation}
    \bu_t = \mathcal{F}([\bu, \bw]) \in \bV^{*} \label{riesz}
\end{equation}
where 
\begin{equation}
\label{F_equation}
\begin{split}
     \langle \mathcal{F}([\bu, \bw]), \bvarphi \rangle &= -(\varepsilon(\bu), \varepsilon(\bvarphi))_{\Omega_f} + \langle \sigma(\bw) \cdot\nu, \bvarphi\rangle_{\Gamma_s},
    \quad \text{for all }\bvarphi \in \bV.
\end{split}
\end{equation}
Secondly, via the fluid PDE component \eqref{eq:fsif1}, we have for all \begin{equation}\bpsi \in \mathcal{V} \equiv \{\bpsi \in \bH_{0}^{1}(\Omega_f) : \div(\bpsi) = 0 \text{ in } \Omega_f \}:\end{equation}
\begin{equation*}
\begin{split}
-(\bu_t - \div(\varepsilon(\bu)), \bpsi)_{\Omega_f} = (\nabla p , \bpsi)_{\Omega_f} = \cancel{\langle p\nu, \bpsi\rangle_{\partial\Omega_f}} - (p, \div(\bpsi))_{\Omega_f} = 0.
\end{split}
\end{equation*}
That is, $\bu_{t} - \div(\varepsilon(\bu)) \in \mathcal{V}^{\perp}$, where
\begin{equation}\mathcal{V}^{\perp} \equiv \{\boldsymbol{\xi} \in \bH_{0}^{-1}(\Omega_f) : \langle \boldsymbol{\xi}, \bpsi\rangle = 0 \text{ for all } \bpsi \in \mathcal{V} \}:\end{equation}
Subsequently, if $p(t) \in L^{2}(\Omega_f)$ has the decomposition $p = q_0 + c_0$, where
\begin{equation}
    q_0 \in \hat{L}^{2}(\Omega_f) \equiv \{ q \in L^{2}(\Omega_f) : \int_{\Omega_f} g \, d\Omega_f = 0 \},
\end{equation}
and
\begin{equation}
    c_0 \equiv \text{constant},
\end{equation}
then from the classic result of de Rham (see e.g. Lemma 2.2.2., p. 75 in \cite{sohr2012navier}), we have the existence of a continuous map, $\mathcal{L} : \mathcal{V}^{\perp} \to \hat{L}^{2}(\Omega_f)$ (where orthogonal complement $\mathcal{V}^{\perp} \subset \bH^{-1}(\Omega_f))$, such that the image

\begin{equation}
\begin{alignedat}{2}
\mathcal{L}(\div(\varepsilon(\bu)) - \bu_t)
  &=\;&&
    q_0\quad\text{in }\hat{L}^{2}(\Omega_f)
\end{alignedat}
\label{orthogonal}
\end{equation}
yields
\begin{equation}
\begin{alignedat}{2}
\bu_t - \div(\varepsilon(\bu)) + \nabla q_0
  &=\;&&
    0\quad\text{in }\Omega_f.
\end{alignedat}
\label{newfluideq}
\end{equation}
Combining \eqref{riesz} and \eqref{orthogonal} then gives
\begin{equation}
\begin{alignedat}{2}
q_0 &= \mathcal{L}[\div(\varepsilon(\bu)) - (\mathcal{F}([\bu, \bw]))]\\
  &= \mathcal{L}(\mathcal{P} \div(\varepsilon(\bu)) - \mathcal{L}(\mathcal{P}  \mathcal{F}([\bu, \bw])).
\end{alignedat}
\end{equation}
(Here, $\mathcal{P} : \bH^{-1}(\Omega_f) \to \mathcal{V}^{\perp}$ is the orthogonal projection onto $\mathcal{V}^{\perp}$.) Note that the second term of the right hand side is well-defined \textit{if} $\mathcal{F}([\bu, \bw]) \in \bL^{2}(\Omega_f)$. If this is so, we can then define the linear maps $\bG_1$ and $\bG_2$ as
\begin{equation}
    \bG_{1}\bu \equiv -\nabla \mathcal{L}(\mathcal{P} \div(\varepsilon(\bu)) ),
\label{G1_map}
\end{equation}
\begin{equation}
    \bG_{2}([\bu, \bw]) \equiv \nabla \mathcal{L} ( \mathcal{P}  \mathcal{F}([\bu, \bw])).
\label{G2_map}
\end{equation}
Consequently these expressions and \eqref{newfluideq} allow writing the fluid PDE component of the system \eqref{eq:fsif1}--\eqref{eq:fsif3} in terms of $\bu$ and $\bw$ alone, thereby eliminating the pressure so as to have:
\begin{equation}
\begin{alignedat}{2}
u_t = \div(\varepsilon(\bu)) + \bG_{1}\bu + \bG_{2}([\bu, \bw])\quad\text{in }(0, T) \times \Omega_f.
\end{alignedat}
\end{equation}

\section{Consequences of new pressure elimination method}
\subsection{The explicit form of the fluid--structure generator}
We write the FSI \eqref{eq:fsif1}--\eqref{eq:fsif8} as the Cauchy problem,

\begin{equation}\frac{d}{dt} \begin{bmatrix}\bu\\ \bw\\ \bw_{t}\end{bmatrix} = \mathcal{A}\begin{bmatrix}\bu\\ \bw\\ \bw_{t}\end{bmatrix},\label{evol1}\end{equation}
\begin{equation}[\bu(0), \bw(0), \bw_{t}(0)] = [\bu_0, \bw_{0}, \bz_{0}] \in \bH,\label{evol2}\end{equation}
where

\begin{equation}\mathcal{A}\begin{bmatrix}\bu\\ \bw\\ \bw_{t}\end{bmatrix} \equiv \begin{bmatrix} \div(\varepsilon(\bu)) + \bG_{1}\bu + \bG_{2}([\bu, \bw]) \\ \bw_t \\ \div\sigma(\bw) - \bw 
\end{bmatrix}\label{operator}.\end{equation}
The domain $D(\mathcal{A}) \subset \bH$ will be specified below, with $\nabla p$ being determined by $[\bu, \bw, \bz] \in D(\mathcal{A})$, with the $\bG_1$ and $\bG_2$ given by \eqref{G1_map} and \eqref{G2_map}.

\bigskip

\noindent Before stating $D(\mathcal{A})$ explicitly, we first prove some preliminaries.

\begin{mylemma}
\label{lemma1}
    Suppose $\bv$ is in $\bH^{1}(\Omega_f)$, $\pi$ is in $L^{2}(\Omega_f)$, and $\div(\varepsilon(\bv)) - \nabla \pi = \mathbf{F} \in \bL^{2}(\Omega_f)$. Then we have that $\varepsilon(\bv) \cdot \nu - \pi \nu \in \bH^{-1/2}(\partial\Omega_f)$, with the following boundary trace estimate:

\begin{equation}\| \varepsilon(\bv) \cdot \nu - \pi \nu |_{\partial\Omega_f} \|_{-1/2, \partial\Omega_f} \, \leq C (\|\bv\|_{1, \Omega_f} + \|\pi\|_{0, 
\Omega_f} + \|\mathbf{F}\|_{0, \Omega_f}).\end{equation}
\end{mylemma}
\noindent \textit{Proof.} Since the (Dirichlet) Sobolev trace map
$$\gamma_0 \in \mathcal{L}(\bH^{1}(\Omega_f), \bH^{1/2}(\partial\Omega_f) )$$
for a general Lipschitz boundary $\partial\Omega$ is surjective, where $\gamma_{0}(\widetilde{\bv}) = \widetilde{\bv}|_{\partial\Omega_f}$ the pointwise restrictions for $\widetilde{\bv} \in \bC^{\infty}(\overline{\Omega_f})$, then there exists a continuous right inverse $\gamma_{0}^{+} \in \mathcal{L}(\bH^{1/2}(\partial\Omega_f), \bH^{1}(\Omega_f) )$. (See e.g., Theorem 3.38, p. 102 of \cite{mclean2000strongly}.) That is, $\gamma_0 \gamma_0^{+}(\widetilde{\bg}) = \widetilde{\bg}$ for all $\widetilde{\bg} \in \bH^{1/2}(\partial\Omega_f)$. Subsequently, we have via an integration by parts:
\begin{equation*} \begin{alignedat}{2} \langle \varepsilon(\bv) \cdot \nu - \pi \nu, \widetilde{\bg} \rangle_{\partial\Omega_f}
&= (\varepsilon(\bv) , \varepsilon(\gamma_{0}^{+}(\widetilde{\bg})) )_{\Omega_f} - (\pi, \div \gamma_{0}^{+}(\widetilde{\bg}))_{\Omega_f} + (\div\varepsilon(\bv) - \nabla \pi, \gamma_{0}^{+}(\bg))_{\Omega_f} \\
&= (\mathbf{F}, \gamma_{0}^{+}(\widetilde{\bg}))_{\Omega_f} + (\varepsilon(\bv), \varepsilon(\gamma_{0}^{+}(\widetilde{\bg})))_{\Omega_f} - (\pi, \div \gamma_{0}^{+}(\widetilde{\bg}))_{\Omega_f}, \end{alignedat} \end{equation*}
where given $\widetilde{\bg} \in \bH^{1/2}(\partial\Omega_f)$. Then:
$$|\langle \varepsilon(\bv) \cdot \nu - \pi \nu, \widetilde{\bg} \rangle_{\partial\Omega_f} | \, \leq C (\|\bv\|_{1, \Omega_f} + \|\pi\|_{0, \Omega_f} + \|\mathbf{F}\|_{0, \Omega_f} ) \|\widetilde{\bg}\|_{1/2, \partial\Omega_f},$$
where we are using the continuity of $\gamma_{0}^{+}$. (Subsequent to this proof, we will hence denote the trace $\gamma_{0}(\bv)$ as $\bv|_{\partial\Omega}$.)
\qed

\subsection{\texorpdfstring{Domain of $\mathcal{A}$}{Domain of A}}
We define the domain of $\mathcal{A} : \bH \to \bH$ to be the subspace $D(\mathcal{A})$ of $\bH$ composed of all $[\bu_0, \bw_0, \bz_0] \in \bH$ which satisfy the following:

\begin{enumerate}[leftmargin=4em, itemindent=1em, label=(D.\arabic*)]
    \item \hypertarget{domain1}{} $[\bu_0, \bw_0, \bz_0] \in (\bV \cap \bH_f) \times \bH^{1}(\Omega_s) \times \bH^{1}(\Omega_s)$, where $\bV$ is given in \eqref{test_space}.
    \item \hypertarget{domain2}{} The structural component $\bw_0$ satisfies $\div(\sigma(\bw_0)) \in \bL^{2}(\Omega_s)$. (Subsequently, a straightforward energy argument yields that $\sigma(\bw_0) \cdot \nu |_{\Gamma_s}$ is well-defined as an element of $\bH^{-1/2}(\Gamma_s)$; see, e.g., p. 115, Lemma 4.3 of \cite{mclean2000strongly}.)
    \item \hypertarget{domain3}{} The components obey the following relation on the boundary interface $\Gamma_s$:
    \begin{equation*} \begin{alignedat}{2} \bu_{0}|_{\Gamma_s} &=\;&& \bz_{0}|_{\Gamma_s}\quad\text{on }\Gamma_s. \end{alignedat} \end{equation*}
    \item \hypertarget{domain4}{} For the given data $[\bu_0, \bw_0, \bz_0]$, there exists a corresponding pressure function $\pi_0 \in L^{2}(\Omega_f)$ such that:
    \begin{enumerate}[leftmargin=4em, label=(D.4.\alph*)]
        \item \hypertarget{domain4a}{} The pair $(\bu_0, \pi_0)$ satisfies:
        \begin{equation} -\div(\varepsilon(\bu_0)) + \nabla \pi_0 \in \bH_f. \label{generator_a}\end{equation}
        Consequently an integration by parts (see Lemma~\ref{lemma1}) yields:
        \begin{equation}\varepsilon(\bu_0)\cdot\nu - \pi_{0}\nu \in \bH^{-1/2}(\partial\Omega_f),\label{generator_a2}\end{equation}
        and so $\varepsilon(\bu_0) \cdot \nu - q_{0}\nu \in \bH^{-1/2}(\partial\Omega_f)$, where \begin{equation}\label{derhamsplit}\pi_0 = q_0 + c_0, \text{ where } \int_{\Omega_f} q_0 \, d\Omega_f = 0, \quad c_0 = \text{constant}.\end{equation}
        \item \hypertarget{domain4b}{} One has the $\Gamma_s$ boundary interface condition on the components $[\bu_0, \bw_0]$ and associated pressure function $\pi_0$
        \begin{equation}\sigma(\bw_0) \cdot \nu = \varepsilon(\bu_0) \cdot \nu - \pi_{0}\nu.\label{generator_b}\end{equation}
        \end{enumerate}
        We note that the interface condition \eqref{derhamsplit} and domain criterion (\hyperlink{domain4a}{D.4a}) yield that \begin{equation}\mathcal{F}([\bu_{0},\bw_{0}])\in \bH_{f},\end{equation} where $\mathcal{F}(\cdot)$ is as given in \eqref{F_equation}. Indeed, in fact from (\hyperlink{domain4a}{D.4a}) there exists $\pi _{0}\in L^{2}(\Omega_{f})$ such that 
        \begin{equation*}
        -\div(\varepsilon (\bu_{0}))+\nabla \pi _{0}=\mathbf{f}\text{, say, in }\bH_{f}. \end{equation*} Multiplying both sides of this equation by $\bvarphi \in \mathbf{V}$, integrating and then integrating by parts, we have 
        \begin{equation*}-\left\langle \varepsilon (\bu_{0})\cdot\nu - \pi_0 \nu,\bvarphi \right\rangle _{\Gamma_s}+\left( \varepsilon (\bu_{0}),\varepsilon (\bvarphi )\right) _{\Omega _{f}}=\left( \mathbf{f},\bvarphi  \right) _{\Omega _{f}}.\end{equation*}
        Using the matching fluxes BC between elastic components in \eqref{generator_b}, this relation then becomes 
        \begin{equation*}
        \left\langle \mathcal{F}([\bu_{0},\bw_{0}]),\bvarphi \right\rangle =\left( \mathbf{f},\bvarphi \right) _{\Omega _{f}} \text{, for all }\bvarphi \in \bV\text{,} \end{equation*}
        where $\mathcal{F}$ is as given in \eqref{F_equation}. Since $\bV$ is dense in $\mathbf{H}_{f}$, we then obtain the conclusion that $\mathcal{F}([\bu_{0},\bw_{0}])\in \mathbf{H}_{f}$, as required. Thus, $q_{0}\in \hat{L}^{2}(\Omega )$ of \eqref{derhamsplit} is necessarily given by
        \begin{equation*} q_{0}=\mathcal{L}(\mathcal{P}\div(\varepsilon (\bu_{0}))-\mathcal{L}(\mathcal{PF}([\bu_{0},\bw_{0}])).\end{equation*}
        Subsequently, the constant component of the pressure in \eqref{derhamsplit} is recovered via the matching fluxes boundary condition in \eqref{generator_b}:\begin{equation*}
        c_{0}=\left[ \varepsilon (\bu_0)\cdot \nu-\sigma (\bw_0)\cdot\nu\right] _{\Gamma _{s}}\cdot \nu-\left.q_{0}\right\vert _{\Gamma s}.\end{equation*}
\end{enumerate}

\subsection{\texorpdfstring{Semigroup well-posedness of \eqref{eq:fsif1}--\eqref{eq:fsif8}}{Semigroup well-posedness}}
\label{wellposedness}
For the fluid--structure model that is defined in \eqref{eq:fsif1}--\eqref{eq:fsif8} and (\hyperlink{domain1}{D.1})--(\hyperlink{domain4}{D.4}), we have semigroup generation for the associated operator $\mathcal{A}$ and domain defined in $D(\mathcal{A})$:

\begin{mytheorem} With respect to bounded Lipschitz pair $\{ \Omega_f, \Omega_s\}$, let modeling operator $\mathcal{A} : \bH \to \bH$ of \eqref{eq:fsif1}--\eqref{eq:fsif8} be as given in \eqref{operator}, with domain $D(\mathcal{A})$ as given in (\hyperlink{domain1}{D.1})--(\hyperlink{domain4}{D.4}). Then:
\begin{enumerate}[label=(\roman*)]
\item The operator $\mathcal{A} : D(\mathcal{A}) \subset \bH \to \bH$ generates a contraction $C_0$-semigroup $\{e^{\mathcal{A}t} \}_{t \geq 0}$ on $\bH$. That is, for $[\bu_0, \bw_0, \bz_0 ] \in \bH$, the solution to $[\bu, \bw, \bw_t ] \in \mathbf{C}([0, T]; \bH)$ of \eqref{eq:fsif1}--\eqref{eq:fsif8} is given by:

$$[\bu(t), \bw(t), \bw_{t}(t)] = e^{\mathcal{A}t} [\bu_0, \bw_0,\bw_1].$$

Furthermore, we have the additional regularity for the fluid component:
\begin{equation} \bu \in L^{2}(0, T; \bH^{1}(\Omega_f)).\label{extra_regularity}\end{equation}

\item If $[\bu_0, \bw_0, \bz_0 ] \in D(\mathcal{A})$, one obtains the additional regularity for the solution $[\bu, \bw, \bw_t]$ of \eqref{eq:fsif1}--\eqref{eq:fsif8}: $[\bu, \bw, \bw_t] \in \mathbf{C}([0, T]; D(\mathcal{A})), p \in C([0, T]; L^{2}(\Omega_f))$, with pressure $p$ given explicitly by:
\begin{equation}
\begin{aligned}
p &= q_0 + c_0, \\
&\quad \text{where } \begin{aligned}[t]
    q_0 &= \mathcal{L}(\mathcal{P} \div(\varepsilon(\bu)) - \mathcal{L}(\mathcal{P} [\mathcal{F}([\bu, \bw])]), \\
    c_0 &\equiv \left[ \varepsilon (\bu)\cdot \nu-\sigma (\bw)\cdot\nu\right] _{\Gamma _{s}\cdot \nu} -\left.q_{0}\right\vert _{\Gamma s}.
\end{aligned}
\end{aligned}
\end{equation}

\end{enumerate}

\end{mytheorem}

The proof of Theorem 6, apart from the new pressure elimination method and Proposition~\ref{newproposition}, closely mirrors that in \cite{avalos2008maximality}, undertaken for smooth pair $\{\Omega_f, \Omega_s\}$. For the sake of completion, we provide the total steps for semigroup generation in Appendix~\ref{appendix:semigroup_proof}).
The proof for the Theorem essentially follows from: completing most of the steps of proof for semigroup generation as in \cite{avalos2008maximality} (with only minor adjustments for the reduced trace regularity, see Appendix~\ref{appendix:semigroup_proof}). Specifically, we show that the fluid--structure generator is maximal dissipative (see Appendix~\ref{appendix:semigroup_proof} sections \ref{appendix:dissipative proof}--\ref{appendix:maximality_proof_2} for verification of dissipativity and maximality), and contraction semigroup generation follows by application of the Lumer--Phillips Theorem. Except for Proposition~\ref{newproposition} below, the proof of maximality follows the same steps as in \cite{avalos2008maximality}. The proof of the extra regularity \eqref{extra_regularity} follows via applying standard energy methods to the system \eqref{eq:fsif1}--\eqref{eq:fsif8}, by means of multiplying the fluid component by $\bu$, the elastic equation by $\bw_{t}$ and integration in time and space.

In the preceding, we note that a key point of departure from the maximal dissipative argument in \cite{avalos2008maximality} is the need to establish the associated ``inf-sup'' estimate on Lipschitz domain $\Omega_f$. The inequality ultimately allows for an application of the Babu\v{s}ka-Brezzi Theorem (see Theorem~\ref{LBB}), by way of showing the maximality of $\mathcal{A}: \bH \to \bH$. To wit, we have (in the style of \cite{dvorak2008_phd}):

\begin{myproposition}[Inf-sup property for general Lipschitz domains]
\label{newproposition}
Define the bilinear form $b(\cdot, \cdot) : \bSigma \times L^{2}(\Omega_f) \to \mathbb{R}$, where $\bSigma \equiv \bH_{\Gamma_f, 0}^{1}(\Omega_f) = \{\bvarphi \in \bH^{1}(\Omega_f) : \bvarphi|_{\Gamma_f} = 0 \}$:
\begin{equation}
\begin{alignedat}{2}
b(\bvarphi, \mu)
  &\equiv -(\mu, \div\bvarphi)_{\Omega_f}
  &\quad& \text{ for all }\,\bvarphi \in \bSigma, \mu \in L^{2}(\Omega_f).
\end{alignedat}
\label{bilinear_infsup}
\end{equation}
Then we have satisfaction of the inf-sup condition \eqref{infsuplemma} for $b(\cdot, \cdot)$: namely there exists a constant $\beta > 0$ such that
\begin{equation}
\label{infsuplemma}
\sup_{\bvarphi \in \bSigma}\frac{b(\bvarphi, \mu)}{\|\bvarphi\|_{\bSigma} } \geq \beta \|\mu\|_{0, \Omega_f} \quad \text{ for every } \mu \in L^{2}(\Omega_f).
\end{equation}
\end{myproposition}

\bigskip

\noindent \textit{Proof.} First we recall the fact that in any open bounded subset $\Omega$ of $\mathbb{R}^n$ with a Lipschitz boundary $\Gamma$, there exist $\delta > 0$ and $\boldsymbol{\zeta} \in \bC^{\infty}(\overline{\Omega})$ such that $\boldsymbol{\zeta} \cdot \nu \geq \delta$ a.e. on $\Gamma$ (see \cite{grisvard1985_elliptic}), Lemma 1.5.1.9). Suppose $\mu \in L^{2}(\Omega_f)$ is given, and let $\bomega \in \bSigma$ be the solution of the BVP:
\begin{equation}
\begin{cases} 
\mathrm{div}(\bomega) = -\mu \langle \boldsymbol{\zeta}, \nu \rangle_{\Gamma_s} & \text{ in } \Omega_f,
\\
\bomega|_{\Gamma_f} = 0 & \text{ on } \Gamma_f,
\\
\bomega|_{\Gamma_s} = -\left(\int_{\Omega_f} \mu\, d\Omega_f\right)\boldsymbol{\zeta}(x) & \text{ on } \Gamma_s.
\end{cases}
\end{equation}
We have existence of this solution for a general Lipschitz domain by p.127 of \cite{galdi1994navierstokes}, using satisfaction of the compatibility conditions, and also have the following estimate for some positive constant $C$:
\begin{equation}\| \nabla \bomega \|_{0, \Omega_f} \; \leq \; C \|\mu\|_{0, \Omega_f}.\end{equation}
We then have, for the given $L^{2}$-function $\mu$, 
\begin{equation}
\begin{aligned}
\sup_{\bvarphi \in \bSigma} \frac{b(\bvarphi, \mu)}{\|\nabla \bvarphi\|_{0, \Omega_f}} 
&\;=\; \sup_{\bvarphi \in X} \frac{-\int \mu \, \mathrm{div}(\bvarphi) \, d\Omega_f}{\|\nabla \bvarphi\|_{0, \Omega_f}} \;\geq\; \frac{-\int \mu \, \mathrm{div}(\bomega) \, d\Omega_f}{\|\nabla \bomega\|_{0, \Omega_f}} \\
&\;=\; \frac{\int \mu^{2} \langle \boldsymbol{\zeta}, \nu \rangle_{\Gamma_s} \, d\Omega_f}{\|\nabla \bomega\|_{0, \Omega_f}} \;\geq\; \delta \cdot \text{meas}(\Gamma_s)\frac{\|\mu\|_{0, \Omega_f}^2}{\|\nabla \bomega\|_{0, \Omega_f}} \\& \;\geq\; \frac{\delta \cdot \text{meas}(\Gamma_s)}{C} \|\mu\|_{0, \Omega_f}.
\end{aligned}
\end{equation}
Thus the inf-sup condition \eqref{infsuplemma} is satisfied (using norm equivalence).
\qed

\bigskip

Subsequently, we apply the remarks in Section~\ref{newmethod} for recovery of the pressure. We emphasize that the proof is enabled by the new definition of both $\mathcal{A}$ in \eqref{operator} and $D(\mathcal{A})$, which was made possible by the new pressure elimination method.

\par Then application of the Sobolev Trace Theorem and the and boundary condition in \eqref{eq:fsif5} gives the following Corollary:

\begin{mycorollary}
Given $[\bu_0, \bw_0, \bz_0 ] \in \bH$, the mechanical velocity component of the solution $[\bu, \bw, \bw_t]$ to \eqref{eq:fsif1}--\eqref{eq:fsif8} has the regularity: \begin{equation}\bw_{t}|_{\Gamma_s} \in L^{2}(0, T; \bH^{1/2}(\Gamma_s)).\end{equation}
\end{mycorollary}

\subsection{Convergence estimates for approximating FEM}
\label{FEM_estimates}

We now state the setting in \cite{avalos2016infsup} that our approximation scheme satisfies. Assume $\Omega \subset \mathbb{R}^d$, $d = 2, 3$, and write $\bY = [\bu, \bw, \bz] \in \bH$ so that the evolution system \eqref{evol1}--\eqref{evol2} is:
\begin{equation}\bY_{t}(t) = \mathcal{A}\bY(t), \qquad \bY(x) = \begin{bmatrix}\bu_0\\ \bw_0\\ \bz_0\end{bmatrix} \in \bH. \end{equation}

\noindent The associated static PDE (resolvent) problem in \eqref{resolvent} is finding $\bY \in D(\mathcal{A})$---this solution variable now taken to be steady-state---such that:
\begin{equation}(\lambda - \mathcal{A})\bY = \bY^{*},\end{equation}
for $\lambda > 0$ and given $\bY^{*} = (\bu^{*}, \bw^{*}, \bz^{*}) \in \bH$. This reduces to finding the velocity and pressure pair $[\bu, \pi] \in \bH_{\Gamma_f, 0}^{1}(\Omega_f) \times L^{2}(\Omega_f)$ that solves:
\begin{equation}
\begin{alignedat}{2}
a_{\lambda}(\bu,\bvarphi) \;+\; b(\bvarphi,\pi)
  &= F(\bvarphi)
  &\quad& \text{ for all }\,\bvarphi\in \bH_{\Gamma_f,0}^{1}(\Omega_f),\\
b(\bu,\mu)
  &= 0
  &\quad& \text{ for all }\mu\in L^{2}(\Omega_f).
\end{alignedat}
\label{mixed_variational_0}
\end{equation}
as derived in \eqref{mixed_variational} of the Appendix.

Now let $\{ \mathcal{T}_h \}_{h > 0}$ denote a quasi-uniform (\cite{ern2004fem}, Def 1.140, p.76]) family of affine meshes on $\Omega_f \cup \overline{\Omega}_s$ such that every element $K \in \mathcal{T}_h$ resides in either $\overline{\Omega}_f$ or in $\overline{\Omega}_s$. Here, as usual, each subscript $h$ denotes the ``level of refinement of the mesh'' $\mathcal{T}_{h}$ (see \cite{ern2004fem}, p. 32). That is, if given $K \in \mathcal{T}_{h}$, $h_{K} \equiv \textrm{diam}(K)$, then \begin{equation}h \equiv \max_{K \in \mathcal{T}_{h}} h_{K}.\end{equation} Relative to $\{ \mathcal{T}_{h}\}$, we consider approximations of $\bSigma \equiv \bH_{\Gamma_f, 0}^{1}(\Omega_f)$, $\bW \equiv \bH_{0}^{1}(\Omega_s)$ and $\Pi\equiv L^{2}(\Omega_f)$ based on triangular or tetrahedral Taylor--Hood elements $\mathbb{P}_2 / \mathbb{P}_1$ (see \cite{ern2004fem}). For $K \in \mathcal{T}_h$, we let $P^{m}(K)$ denote the set of polynomials of degree $\leq m$ on $K$. Then define approximating spaces by:
\begin{align}
    \bSigma_h \equiv \{ \bu_h \in \mathbf{C}^{0}(\Omega_f) : \bu_{h}|_{\Gamma_f} = 0 \text{ and for all } K \in \mathcal{T}_h \cap \overline{\Omega}_f, \, \bu_{h}|_{K} \in P^{2}(K) \}; \\
    \bW_h \equiv \{ \bw_h \in \mathbf{C}^{0}(\Omega_s) : \text{ for all } K \in \mathcal{T}_h \cap \overline{\Omega}_s , \, \bw_{h}|_{K} \in P^{2}(K) \}; \\
    \Pi_h \equiv \{ \mu_h \in C^{0}(\Omega_f) : \text{ for all } K \in \mathcal{T}_h \cap \overline{\Omega}_f,  \, \mu_{h}|_{K} \in P^{1}(K) \}.
\end{align}

\paragraph{Semi-discrete scheme.} This refers to the discretization of the Hilbert spaces $\bSigma$, $\bW$ and $\Pi$ solely, not the bilinear forms. To wit, the approximate solution is the pair $(\bu_h, \pi_h) \in \bSigma_n \times \Pi_{h}$ that uniquely solves the mixed variational system:
\begin{equation}
\begin{alignedat}{2}
a_{\lambda}(\bu_h, \bvarphi_h) + b(\bvarphi_h, \pi_h)
  &= F(\bvarphi_h)
  &\quad& \text{ for all }\,  \bvarphi_h \in \bSigma_h,
  \\b(\bu_h, \mu_h)
  &= 0
  &\quad& \text{ for all }\,  \mu_h \in \Pi_h,
\end{alignedat}
\label{semi_discrete_scheme_ref}
\end{equation}
where bilinear forms $a_{\lambda}(\cdot, \cdot)$, $b(\cdot)$ are as in \eqref{bilinear}--\eqref{bilinear2} of the Appendix. With respect to $\{\bSigma_h, \Pi_h \}$ and discrete system \eqref{semi_discrete_scheme_ref} we present the following critical result for $\mathbb{R}^2$ and triangular elements with $\Omega_f \subset \mathbb{R}^2$ and each $K \in \mathcal{T}_h$ a triangle, although as noted in \cite{avalos2016infsup}, there are no obstacles to its extension to three dimensions.

\begin{mytheorem}[Discrete uniform inf-sup inequality,  see \cite{avalos2016infsup} Theorem 3.1]
\label{discreteinfsup}
Assume that for each $h > 0$, every element in $\mathcal{T}_h$ that is supported in $\overline{\Omega}_f$ has at least one vertex not in $\Gamma_f$. Then there is a constant $C^{*} > 0$ independent of $h \in (0, h_0 )$, some $h_0 > 0$, such that
\begin{equation}S_{\mu_h} \equiv \sup_{0 \neq \bvarphi_h \in \bSigma_h} \frac{b(\bvarphi_h, \mu_h)}{|\bvarphi_h|_{1, \Omega_f} } \geq C^{*},\end{equation}
for every $\mu_h \in \Pi_h$, with $\|\mu_h \|_{0, \Omega_f} = 1$. (In this section the $|\cdot|_{k, \Omega}$ notation denotes the Sobolev semi-norm over $\Omega$.) Equivalently, via the Poincar\'e inequality, there is $C_{\pi}^{*} > 0$ such that for all $h \in (0, h_0 )$,
\begin{equation}
\inf_{0 \neq \mu_h \in \Pi_h}\ \sup_{0 \neq \bvarphi_h \in \bSigma_h}\ 
\frac{b(\bvarphi_h,\mu_h)}{\|\bvarphi_h\|_{1, \Omega_f}\,\|\mu_h\|_{0, \Omega_f}}
 \ge  C_{\pi}^{*}.
\end{equation}
\end{mytheorem}
The availability of the uniform discrete inf-sup constant $C^{*}$ for the approximating system \eqref{semi_discrete_scheme_ref} provides, in turn, the sort of FEM convergence estimates which are well-known to hold with respect to the FEM for elliptic problems. To wit, we have the following (see e.g., Theorem 2.34, p. 100 of \cite{ern2004fem} and Corollary 1 of \cite{avalos2016infsup}):
\begin{mycorollary}
If functions $\bu \in \bSigma$ and $\pi \in \Pi$ solve \eqref{mixed_variational_0} and $\bu_h \in \bSigma_h$, $\pi_h \in \Pi_h$ solve \eqref{semi_discrete_scheme_ref}, then under the assumptions of Theorem~\ref{discreteinfsup}, there exists $c > 0$ independent of $h > 0$ small, such that:
\begin{equation}\|\bu - \bu_h \|_{1, \Omega_f} \, + \, \|\pi - \pi_h\|_{0, \Omega_f} \;\leq \; c \bigg( \inf_{\bvarphi_h \in \bSigma_h} \|\bu - \bvarphi_h \|_{1, \Omega_f} + \inf_{\mu_h \in \Pi_h} \|\pi - \mu_h \|_{0, \Omega_f} \bigg).\end{equation}
\end{mycorollary}

\newpage
\section{Numerical results}
\label{numerical_results}
\subsection{Numerical approximation framework for the static FSI}

Here, we consider a finite-dimensional approximation $[\bu_h, \bw_h, \bz_h]$ to the solution $[\bu, \bw, \bz]$ of the resolvent equation \eqref{resolvent} below, for given $[\bu^{*}, \bw^{*}, \bz^{*}]$ of $\bH$, where $h$ is the parameter of discretization with respect to the geometry ${\Omega_f, \Omega_s}$. In this Ritz--Galerkin framework (see for example \cite{axelsson1984finite}):
\begin{enumerate}[label=(\roman*)]
    \item Let $\bSigma_h \subset \bH_{\Gamma_f, 0}^{1}(\Omega_f)$ be a finite-dimensional approximating subspace of $\bSigma = \bH_{\Gamma_f, 0}^{1}(\Omega_f)$;
    \item Let $\Pi_{h} \subset L^{2}(\Omega_f)$ be an approximating subspace of $\Pi = L^{2}(\Omega_f)$;
    \item Let $\bW_h \subset \bH_{0}^{1}(\Omega_s)$ be an approximating subspace of $\bW = \bH_{0}^{1}(\Omega_s)$.
\end{enumerate}
(Eventually, we will specify $\{\bSigma_h, \Pi_h, \bW_h \}$ to be the Taylor--Hood scheme of Section \ref{FEM_estimates}, relative to mesh family $\{\mathcal{T}_h \}_{h > 0}$ .)

By way of numerically approximating the solution to \eqref{resolvent}, then the principal task is to find  $(\bu_h, \pi_h) \in \bSigma_h \times \Pi_h$ which uniquely solves the mixed variational form:
\begin{equation}
\begin{alignedat}{2}
a_{\lambda}(\bu_h, \bvarphi_h) \,+\, b(\bvarphi_h, \pi_h)
  &= F(\bvarphi_h)
  &\quad& \text{ for all }\,  \bvarphi_h \in \bSigma_h,
  \\b(\bu_h, \mu_h)
  &= 0
  &\quad& \text{ for all }\,  \mu_h \in \Pi_h,
\end{alignedat}
\label{semidiscrete}
\end{equation}
\noindent where the bilinear forms $a_{\lambda}( \cdot , \cdot ) : \bH_{\Gamma_f, 0}^{1}(\Omega_f) \times \bH_{\Gamma_f, 0}^{1}(\Omega_f) \to \mathbb{R}$ and $b( \cdot, \cdot) : \bH_{\Gamma_f, 0}^{1}(\Omega_f) \times L^{2}(\Omega_f) \to \mathbb{R}$ are defined as in \eqref{bilinear}--\eqref{bilinear2}. Note that we make no assumptions that $\bSigma_h$ consists of divergence free functions, which is a key advantage of such a mixed variational formulation (see \cite{brezzi1991mixed}).

We can subsequently recover the approximation $\bw_h$ to the structural PDE solution component $\bw$ of \eqref{resolvent} via the resolvent relations \eqref{structure_bvp}. Namely, we find unique $\bw_h \in \bW_h$ which solves
\begin{equation}
\label{structure_bvp_recovery}
\begin{cases}
\begin{aligned}
\lambda^{2}\,(\bw_h, \bpsi_h)_{\Omega_s} 
    &+ (\sigma(\bw_h), \varepsilon(\bpsi_h))_{\Omega_s} 
    + (\bw_h, \bpsi_h)_{\Omega_s} = (\lambda\,\bw^* + \bz^*, \bpsi_h)_{\Omega_s}, \\
    & \text{for all } \bpsi \in \bW_h,
\end{aligned} \\[6pt]
\bw_{h}|_{\Gamma_s} 
    = \frac{1}{\lambda}\,\bigl(\bu_h + \bw^*\bigr)|_{\Gamma_s},
    \quad \text{on } \Gamma_s.
\end{cases}
\end{equation}In turn, we obtain an approximation $\bz_h$ to the structural term $\bz$ of \eqref{resolvent} by means of the relation:
\begin{equation}\bz_h = \lambda \bw_h - \bw^{*}.\end{equation}

\bigskip

Remark: we note that although we have focused here on the numerical resolution of the static resolvent equation \eqref{resolvent}, in principle we can consider numerically approximating \eqref{eq:fsif1}--\eqref{eq:fsif8} via finite difference time discretizations, applied to \eqref{semidiscrete}, these differences predicated on the classic exponential formula (see e.g. \cite{pazy1983semigroups}):

$$\begin{bmatrix}\bu(t)\\ \bw(t)\\ \bw_{t}(t)\end{bmatrix} = e^{\mathcal{A}t}\begin{bmatrix}\bu_0\\ \bw_0\\ \bz_0\end{bmatrix} = \lim_{n \to \infty}\bigg(I - \frac{t}{n}\mathcal{A}\bigg)^{-n} \begin{bmatrix}\bu_0\\ \bw_0\\ \bz_0\end{bmatrix} \quad   \text{ for} \begin{bmatrix}\bu_0\\ \bw_0\\ \bz_0\end{bmatrix} \in \bH. $$

\newpage

\subsection{A numerical example}

\begin{figure}[h]
  \centering
  \includegraphics[scale=0.9]{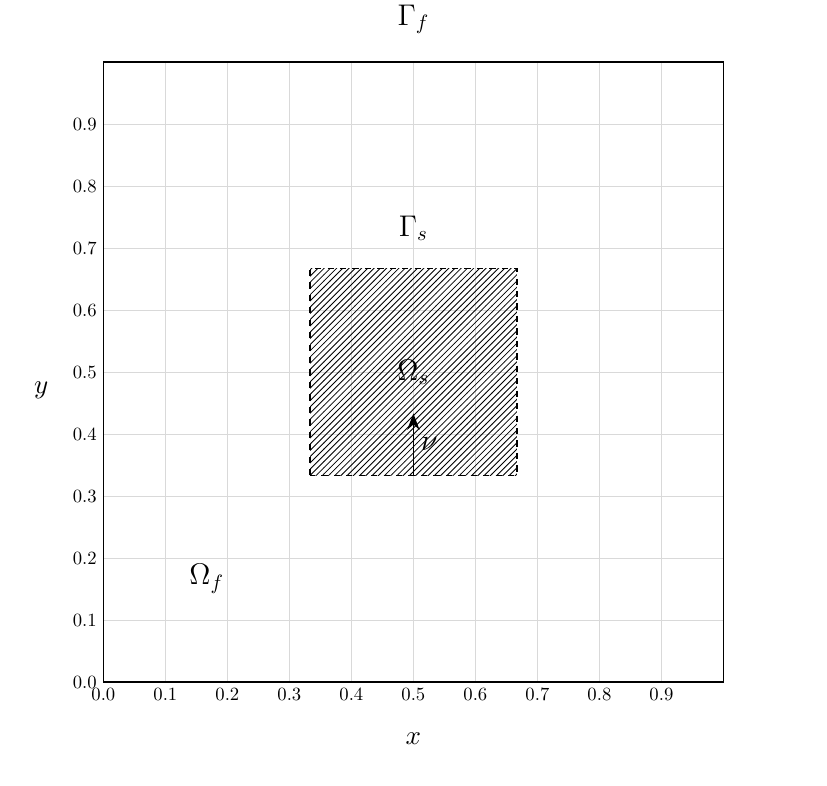}
  \vspace{-0.5cm}
  \caption{The fluid--structure geometry.}
  \label{fig:fsi_geometry_2}
\end{figure}

We now consider a specific example where the solid domain is $\Omega_s = (1/3, 2/3)^2$ and fluid domain $\Omega_f = (0, 1)^2 \setminus [1/3, 2/3]^2$, see Figure~\ref{fig:fsi_geometry_2}.

We will consider the task of approximating the solution of the linearized mixed variational form. First we define the data using the scalar polynomial functions
\begin{align}
A(x) &= \phi(x) \equiv x^{2}(1 - x)^{2}
       \left(x - \tfrac{1}{3}\right)^{3}
       \left(\tfrac{2}{3} - x\right)^{3}, \\
B(y) &= \phi(y),
\end{align}
\noindent where the data we use is:
\begin{align}
(\bu^{*}, \bw^{*}, \bz^{*}) = (\boldsymbol{\rho}, \mathbf{0}, \mathbf{0}) \text{ ; }
\end{align}
where each component of the fluid data $\boldsymbol{\rho}(x, y)$ is:
\begin{align}
\rho_{1}(x, y) &= \lambda A(x) B'(y) - \frac{1}{2}(A''(x)B'(y) + A(x)B'''(y)), \\
\rho_{2}(x, y) &= -\lambda A'(x)B(y) + \frac{1}{2}(A'''(x)B(y) + A'(x)B''(y)),
\end{align}
and where the derivatives are:
\bigskip 
\begin{flalign}
\qquad \phi'(x) &= -10x^{9} + 45x^{8}
       - \tfrac{256}{3}x^{7} + \tfrac{266}{3}x^{6}
       - \tfrac{494}{9}x^{5} + \tfrac{185}{9}x^{4} 
       - \tfrac{3272}{729}x^{3} &\\ \notag
      &\quad + \tfrac{124}{243}x^{2}
       - \tfrac{16}{729}x, &\\
\qquad \phi''(x) &= -90x^{8} + 360x^{7}
        - \tfrac{1792}{3}x^{6} + 532x^{5}
        - \tfrac{2470}{9}x^{4} + \tfrac{740}{9}x^{3} 
        - \tfrac{3272}{243}x^{2} &\\ \notag
       &\quad + \tfrac{248}{243}x
        - \tfrac{16}{729}, &\\
\qquad \phi'''(x) &= -720x^{7} + 2520x^{6}
         - 3584x^{5} + 2660x^{4}
         - \tfrac{9880}{9}x^{3} + \tfrac{740}{3}x^{2} 
         - \tfrac{6544}{243}x &\\ \notag
        &\quad + \tfrac{248}{243}. &
\end{flalign}

\noindent This corresponds to the exact solution of the mixed variational form constructed by the stream function $\psi(x, y)$, and with the choice of pressure $\pi = 0$:
\begin{flalign*}
\qquad  \psi(x, y) &= A(x)B(y), &\\
\bu &= (u_1, u_2) = (\psi_y, -\psi_x ) = (A(x)B'(y), -A'(x)B(y) ).
\end{flalign*}

\begin{figure}[h]
    \centering
    \includegraphics[width=0.7\textwidth]{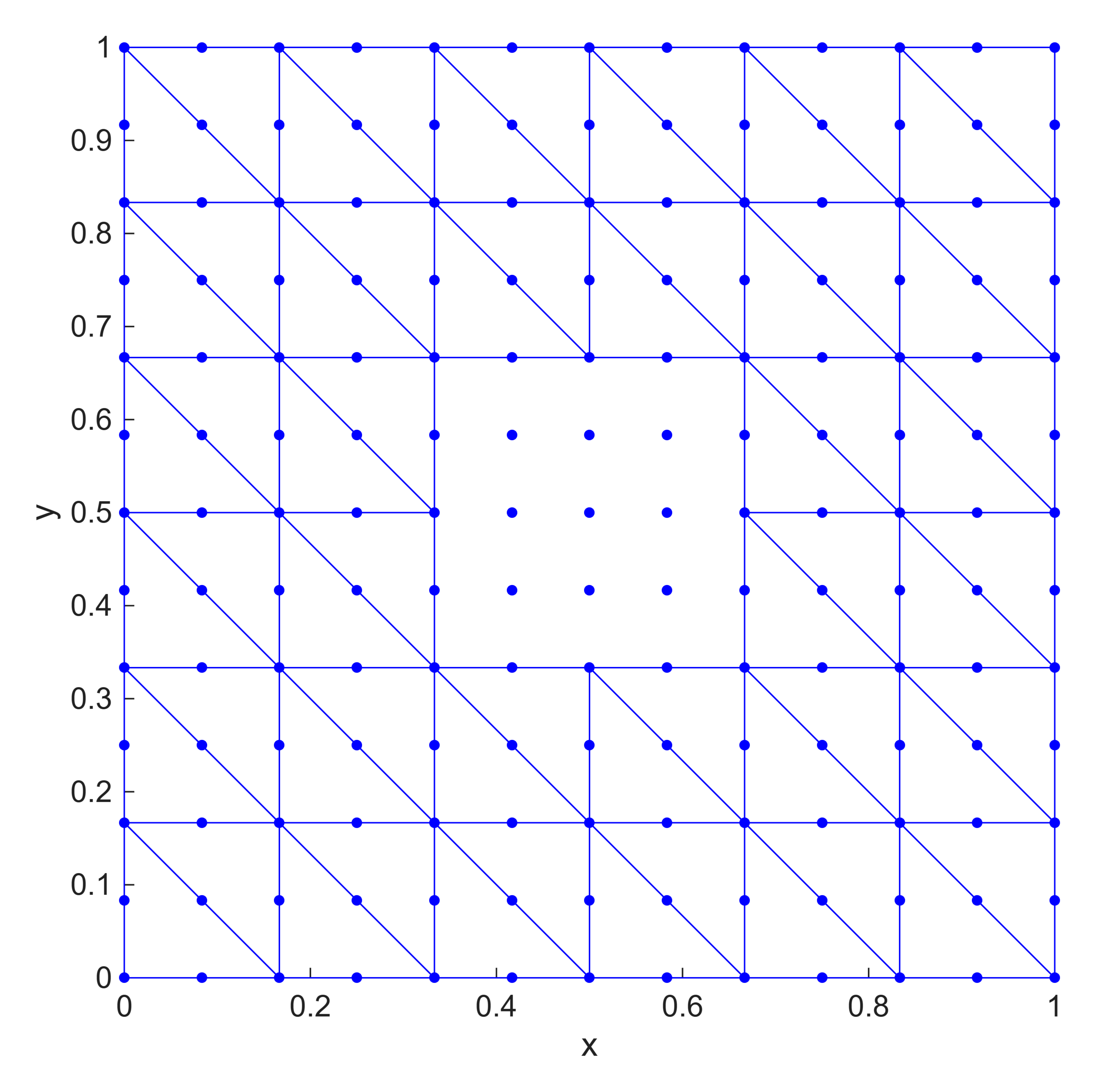}
    \caption{Mesh configuration for the coarsest mesh.}
    \label{fig:mesh1}
\end{figure}

Also note that this corresponds to the exact solution of the static PDE with $(\bu, \bw, \bz)$ with $\bu$ as above, and $\bw = \bz = 0$. We note that the approximation to this example uses non-trivial boundary terms in $a_{\lambda}$ in \eqref{semidiscrete}, thus implicitly making use of the FSI coupling relations \eqref{eq:bc1}--\eqref{eq:bc4} of the static PDE \eqref{resolvent}, with $(\bu, \bw, \bz) \in D(\mathcal{A})$. These latter properties use the vanishing of the higher derivatives of $A(x)$ and $B(y)$ in corresponding regions of the interface $\Gamma_s$.

We use a basic FEM for discretization of the linearized mixed variational form with $\mathbb{P}_2$/$\mathbb{P}_1$ Taylor--Hood elements: the spaces $\bSigma_h$ and $\bW_h$ are the linear spaces of piecewise quadratic basis functions defined on triangular elements, and $\Pi_h$ the linear span of piecewise linear basis functions defined over the same triangular elements. For this test, we choose $\lambda = 1$. We have in Table~\ref{tab:finite-element-error} the errors in the fluid and solid variables for each mesh refinement. (See Figure~\ref{fig:mesh1} for mesh configuration.)

\begin{table}[h]
    \centering
    \caption{Finite element approximation error.} 
    \label{tab:finite-element-error}
    \vspace{0.5em} % vertical space above the table
    \begin{tabular}{ccccc}
        \toprule
        No. of elements & Hypotenuse length & $\|\bu_h - \bu\|_{1, \Omega_f}$ & $\|\pi_h - \pi\|_{0, \Omega_f}$ & $\|\bw_h - \bw \|_{1, \Omega_s}$ \\
        \midrule
        $72$ & $0.235702$ & $5.855 \times 10^{-8}$ & $4.531 \times 10^{-9}$ & $4.617 \times 10^{-10}$ \\
        $288$ & $0.117851$ & $2.965 \times 10^{-8}$ & $5.998 \times 10^{-9}$ & $6.120 \times 10^{-11}$ \\
        $1152$ & $0.0589256$ & $9.342 \times 10^{-9}$ & $1.331 \times 10^{-9}$ & $4.806 \times 10^{-12}$ \\
        $4608$ & $0.0294628$ & $2.488 \times 10^{-9}$ & $1.579 \times 10^{-10}$ & $3.284 \times 10^{-13}$ \\
        $18432$ & $0.00147314$ & $6.331 \times 10^{-10}$ & $1.537 \times 10^{-11}$ & $2.252 \times 10^{-14}$ \\
        $73728$ & $0.00736570$ & $1.590 \times 10^{-10}$ & $1.391 \times 10^{-12}$ & $1.931 \times 10^{-15}$ \\
        \bottomrule
    \end{tabular}
\end{table}

\begin{figure}[H]
    \centering
    \includegraphics[width=0.9\textwidth]{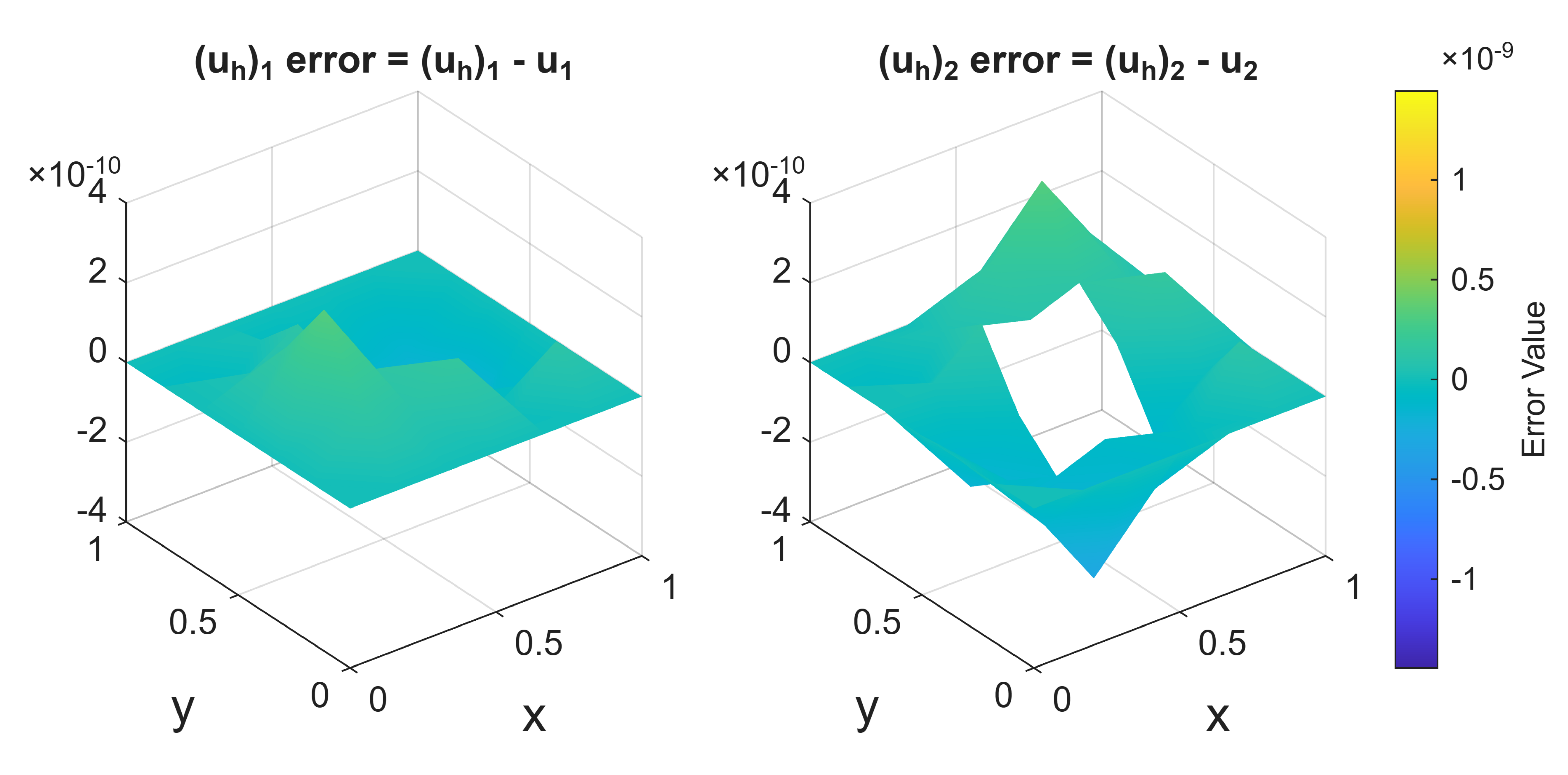}
    \caption{Error for both components of $\bu_h$ compared to $\bu$, for first mesh.}
    \label{fig:mesh4}
\end{figure}

\begin{figure}[H]
    \centering
    \includegraphics[width=0.9\textwidth]{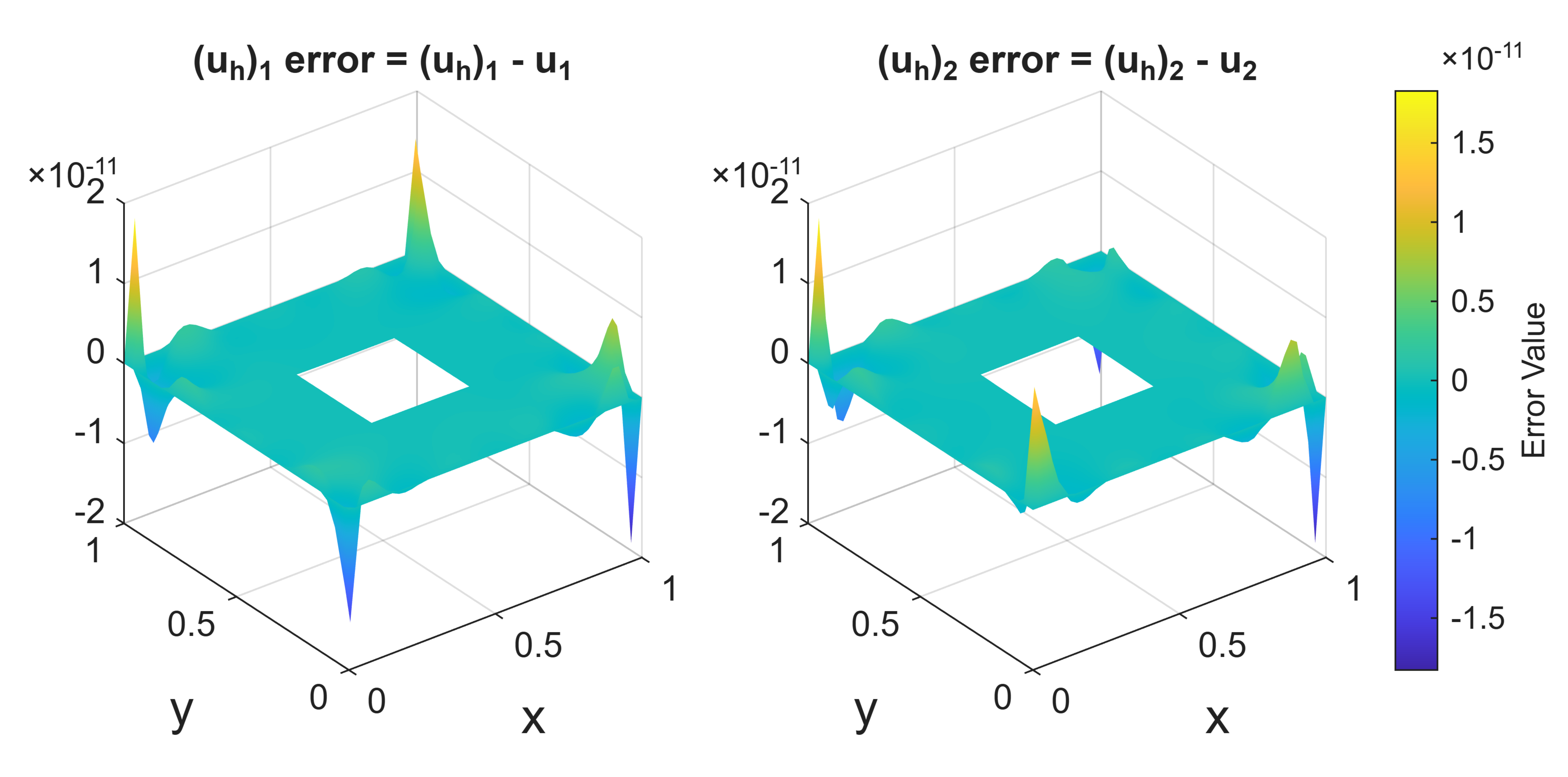}
    \caption{Error for both components of $u_h$ compared to $\bu$, for fourth mesh.}
    \label{fig:mesh5}
\end{figure}

In Figure~\ref{fig:mesh4} and Figure~\ref{fig:mesh5}, we can see that already by the 1st (coarsest) mesh resolution we have the pointwise error within $10^{-10}$, and within $10^{-11}$ by the 4th mesh.

The mesh is refined by a factor of 2 at each step, and so the usual log ratio is simplified into $\log\left(\frac{\text{Error}_{i}}{\text{Error}_{i+1} }\right) / \log(2)$. In the limit this ratio should approach the exponent of convergence $\mathcal{O}(h^k)$. For smooth data and $\mathbb{P}_2/\mathbb{P}_1$ implementation, the expected convergence rate is $k = 2$ for the $\bH^1$ norm of $\bu_h$, and similarly for $\bw_h$. We obtain these rates for the $\bH^1$ norm of $\bu_h$, while the $\bw_h$ and pressure seem to be displaying a better than expected convergence rate that is possibly due to the (unavoidable) approximation (see Table~\ref{tab:error-rates}).

\begin{table}[h]
    \caption{Computed index $k$ in $\mathcal{O}(h^k)$ for FEM.} 
    \centering
    \label{tab:error-rates}
    \vspace{0.5em}
    \begin{tabular}{cccc}
        \toprule
        Meshes & $\bH^{1}$(fluid) & $L^{2}$(pressure) & $\bH^{1}$(solid) \\
        \midrule
        Mesh 1 / Mesh 2 & 0.981 & -0.405 & 2.915 \\
        Mesh 2 / Mesh 3 & 1.666 & 2.172 & 3.671 \\
        Mesh 3 / Mesh 4 & 1.909 & 3.076 & 3.871 \\
        Mesh 4 / Mesh 5 & 1.975 & 3.361 & 3.866 \\
        Mesh 5 / Mesh 6 & 1.993 & 3.466 & 3.544 \\
        \bottomrule
    \end{tabular}
\end{table}

\appendix
\newcounter{appendixchapter}
\setcounter{appendixchapter}{5}
\renewcommand{\theequation}{\theappendixchapter.\arabic{equation}}

\section{Appendix: Proof of semigroup well-posedness}
\label{appendix:semigroup_proof}
\subsection{Proof of dissipativity}
\label{appendix:dissipative proof}

The dissipativity argument follows the same steps as in \cite{avalos2008maximality}. Suppose $[\bu, \bw, \bz] \in D(\mathcal{A})$, then there exists a pressure function $\pi \in L^{2}(\Omega_f)$ so that $[\bu, \bw, \bz, \pi]$ has the properties (\hyperlink{domain1}{D.1})--(\hyperlink{domain4}{D.4}), which lead to the computations below. By definition of $\mathcal{A} : D(\mathcal{A}) \subset \bH \to \bH$ \eqref{operator}, we have:
\begin{equation}
\begin{alignedat}{2}
&\left(\mathcal{A}\begin{bmatrix} \bu \\ \bw \\ \bz\end{bmatrix}, \begin{bmatrix} \bu \\ \bw \\ \bz\end{bmatrix} \right)_{H} =  \left(\begin{bmatrix}\div\varepsilon(\bu) - \nabla \pi\\ \bz\\ \ \div\sigma(\bw) - \bw\end{bmatrix}, \begin{bmatrix} \bu \\ \bw \\ \bz\end{bmatrix} \right)_{H}
\\
 &\quad= (\div\varepsilon(\bu) - \nabla \pi, \bu)_{\Omega_f} + (\varepsilon(\bz), \sigma(\bw))_{\Omega_s} + (\div \sigma(\bw), \bz)_{\Omega_s}
\\
 &\quad= (\div\varepsilon(\bu), \bu)_{\Omega_f} - (\nabla \pi, \bu)_{\Omega_f} + (\div\sigma(\bw), \bz)_{\Omega_s} + (\varepsilon(\bz), \sigma(\bw))_{\Omega_s}.
 \end{alignedat}
 \label{dissip}
\end{equation}

We can deal with each term on the right-hand side:
 \begin{enumerate}[label=(\roman*)]
\item For the first term in \eqref{dissip}, we use Green's theorem and $\bu|_{\Gamma_f} = 0$ from (\hyperlink{domain1}{D.1}) which gives:
\begin{equation}(\div\varepsilon(\bu), \bu)_{\Omega_f} = -\int_{\Omega_f} |\varepsilon(\bu)|^2\, d\Omega_f + \langle \bu, \varepsilon(\bu) \cdot \nu\rangle_{\Gamma_s}\label{dissip1},\end{equation}
where $\langle \cdot , \cdot \rangle$ is the duality pairing of $\bH^{1/2}(\Gamma_s)$ with $\bH^{-1/2}(\Gamma_s)$.

\item Similarly we use Green, (\hyperlink{domain1}{D.1}), the divergence-freeness of (\hyperlink{domain1}{D.1}) and (\hyperlink{domain3}{D.3}) and get:
\begin{equation}-(\nabla \pi, \bu)_{\Omega_f} = -\langle \pi, \bu \cdot \nu\rangle_{\partial\Omega_f} + (\pi, \div(\bu))_{\Omega_f} = -\langle \pi, \bz \cdot \nu \rangle_{\Gamma_s}\label{dissip2}.\end{equation}

\item Using Green's theorem gives
\begin{equation}\label{dissip3}(\div\sigma(\bw), \bz)_{\Omega_s} = -(\sigma(\bw, \varepsilon(\bz))_{\Omega_s} - \langle \bz, \sigma(\bw) \cdot \nu\rangle_{\Gamma_s},\end{equation}
where we use the convention that the normal vector $\nu(x)$ is interior with respect to $\Omega_s$.

\item For the two structure terms, we can use $\bu|_{\Gamma_s} \,= \bz|_{\Gamma_s}$ from (\hyperlink{domain3}{D.3}) and the boundary coupling (\hyperlink{domain4b}{D.4b}) and also \eqref{dissip3} to get the cancellation:
\begin{equation}
\begin{alignedat}{2}
&(\varepsilon(\bz), \sigma(\bw))_{\Omega_s} + (\div\sigma(\bw), \bz)_{\Omega_s}
\\
&\quad = \langle \pi, \bz \cdot \nu\rangle_{\Gamma_s} - \langle \bu, \varepsilon(\bu) \cdot \nu \rangle_{\Gamma_s}.
\end{alignedat}
\label{dissip4}
\end{equation}

\end{enumerate}

Then apply \eqref{dissip1}--\eqref{dissip4} to the right-hand side of \eqref{dissip} to get:
\begin{equation}\left(\mathcal{A}\begin{bmatrix} \bu \\ \bw \\ \bz\end{bmatrix}, \begin{bmatrix} \bu \\ \bw \\ \bz\end{bmatrix} \right)_{H} = -\int_{\Omega_f}|\varepsilon(\bu)|^2 \, d\Omega_f \leq 0,\end{equation}
showing dissipativity of $\mathcal{A}$.

\subsection{Proof of maximality: Mixed variational form}
\label{appendix:maximality_proof_1}
The desired claim is that for sufficiently large $\lambda > 0$, $\Range(\lambda I - \mathcal{A}) = \bH$. Take any $[\bu^{*}, \bw^{*}, \bz^{*}] \in \bH$, a solution $[\bu, \bw, \bz] \in D(\mathcal{A})$ of the equation
\begin{equation}\label{resolvent}(\lambda I - \mathcal{A})\begin{bmatrix}\bu\\ \bw\\ \bz\end{bmatrix} = \begin{bmatrix}\bu^{*}\\ \bw^{*}\\ \bz^{*}\end{bmatrix},\end{equation}
where $\mathcal{A}$ is as defined previously. Componentwise we have:
\begingroup
  \setlength{\jot}{2pt}
  \begin{subequations}
  \begin{align}
    [\,\bu,\bw,\bz\,]\;\in\;\Bigl(\bH^1(\Omega_f) \cap \bH_f\Bigr)
        &\times \bH^1(\Omega_s) \times \bH^1(\Omega_s),
      \label{eq:uvwf}\\
    \lambda\,\bu - \div(\varepsilon(\bu)) + \nabla \pi = \bu^*
      &\quad \text{in }\bH_f,
      \label{eq:fluideq}\\
    \lambda\,\bw - \bz = \bw^*
      &\quad \text{in }\bH^1(\Omega_s),
      \label{eq:wz}\\
    \lambda\,\bz - \div(\sigma(\bw)) + \bw = \bz^*
      &\quad \text{in }\bL^2(\Omega_s).
      \label{eq:zz}
  \end{align}
  \end{subequations}
\endgroup
Then since $[\bu, \bw, \bz] \in D(\mathcal{A})$, we have the additional relations:
\begingroup
  \setlength{\jot}{2pt}
  \begin{subequations}
  \begin{align}
    \bu\bigl|_{\Gamma_f} = 0
      &\quad \text{on }\Gamma_f,
      \label{eq:bc1}\\
    \bu\bigl|_{\Gamma_s} = \bz
      &\quad \text{on }\Gamma_s,
      \label{eq:bc2}\\
    \sigma(\bw)\cdot\nu = \varepsilon(\bu)\cdot\nu - \pi\,\nu
      &\quad \text{on }\Gamma_s,
      \label{eq:bc3}\\
    \div(\bu) = 0
      &\quad \text{a.e.\ in }\Omega_f,
      \label{eq:bc4}
  \end{align}
  \end{subequations}
\endgroup
where the equality for the boundary coupling is taken to be in $\bH^{-1/2}(\Gamma_s)$. We multiply the first structure equation by $\lambda$ to get $\lambda^{2}\bw - \lambda \bz = \lambda \bw^{*}$, and then substitution into the second structure equation gets $\lambda^{2}\bw - \div\sigma(\bw) + \bw = \lambda \bw^{*} + \bz^{*}$. Then notice that $\bw_{t} = \bu \in \Gamma_s$ or equivalently $\bz = \bu|_{\Gamma_s}$, and substitute the first structure equation into this to get $\lambda \bw - \bw^{*} = \bu|_{\Gamma_s}$. Then in total we have the system:

\begin{equation}
\label{structure_bvp}
\begin{cases}
\displaystyle
\lambda^{2}\,\bw -\div(\sigma(\bw)) \, +\,  \bw
  =\lambda\,\bw^* + \bz^*
  & \text{in }\Omega_s,\\[1ex]
\displaystyle
\bw =\frac{1}{\lambda}\,\bigl(\bu + \bw^*\bigr)
  & \text{on }\Gamma_s.
\end{cases}
\end{equation}
We can then write out the solution $\bw$ to this system in terms of the fluid term $\bu$. First define $D_{\lambda} : \bL^{2}(\Gamma_s) \to \bL^{2}(\Omega_s)$ as follows: $\mathbf{g} = D_{\lambda}(\mathbf{h})$ if and only if $\mathbf{g}$ solves the homogeneous elliptic problem:

\begin{equation}
\label{solid_BVP_1}
\begin{cases}
\displaystyle
\lambda^{2} \bw - \div\sigma(\bw)\, + \; \bw
  =0
  & \text{in }\Omega_s,\\[1ex]
\displaystyle
\bw|_{\Gamma_s} = \mathbf{h}
  & \text{on }\Gamma_s,
\end{cases}
\end{equation}
and by elliptic theory, $D_{\lambda} \in \mathcal{L}(\bH^{1/2}(\Gamma_s), \bH^{1}(\Omega_s))$ for $\Omega_s$ a Lipschitz domain. Secondly, we define $\mathcal{L}_{\lambda} : D(\mathcal{L}_{\lambda}) \to \bL^{2}(\Omega_s)$ by:
\begin{equation}\label{solid_BVP_2}\mathcal{L}_{\lambda}\bw \equiv \lambda^{2}\bw - \div\sigma(\bw)\;\, + \; \bw,\end{equation}
with $D(\mathcal{L}_{\lambda}) = \{\bw \in \bH_{0}^{1}(\Omega_s) : \div(\sigma(\bw)) \in \bL^{2}(\Omega_s)\}$. Then by elliptic theory we have $\mathcal{L}_{\lambda}$ is positive self-adjoint, with bounded inverse on $\bL^{2}(\Omega_s)$. We can then re-write the solution by means of the preceding operators:

\begin{equation}\label{structure_via_fluid}
\bw = \frac{1}{\lambda}D_{\lambda}(\bu|_{\Gamma_s}) + \frac{1}{\lambda}D_{\lambda}(\bw^{*}|_{\Gamma_s}) + \mathcal{L}_{\lambda}^{-1}(\lambda \bw^{*} + \bz^{*}) = \bw_1 + \bw_0,\end{equation}
where $\bw_1 = \frac{1}{\lambda}D_{\lambda}(\bu|_{\Gamma_s})$ and the other `data' terms are in $\bw_0$.

\bigskip

Now having expressed $\bw$ in terms of $\bu$ and data terms, note that $\bu$ is still unknown. Now we return to the fluid component above, which solves the system:
\begin{equation*}
\begin{alignedat}{2}
\lambda \bu - \div\varepsilon(\bu) + \nabla p
  &\;=\;&&
    \bu^*\quad\text{in }\Omega_f,
\end{alignedat}
\end{equation*}
and recall $\bu$ vanishes on $\Gamma_f$ and is divergence-free. First define the space of fluid test functions $\bH_{\Gamma_f, 0}^{1}(\Omega_f) = \{\bvarphi \in \bH^{1}(\Omega_f) : \bvarphi|_{\Gamma_f} = 0\} $. We multiply the fluid system above by the test function $\bvarphi \in \bH_{\Gamma_f, 0}^{1}(\Omega_f)$, and use Green's formula for the manipulations below. The idea is that we first compute formally (assuming sufficient trace regularity) to get the mixed variational form \eqref{mixed_variational}, and after applying Theorem~\ref{LBB} (Babu\v{s}ka--Brezzi) gives sufficient regularity of the component terms $\bu$ and $\pi$ which eventually suffices to get maximality. We have:
\begin{equation}(\lambda \bu, \bvarphi) - (\div \varepsilon(\bu), \bvarphi) + (\nabla \pi, \bvarphi) = (\bu^{*}, \bvarphi),\end{equation}
and use of Green's theorem on the term $(\div\varepsilon(\bu), \bvarphi)$ gives:
\begin{equation}
\begin{alignedat}{2}
(\lambda \bu, \bvarphi) + (\varepsilon(\bu), \varepsilon(\bvarphi))_{\Omega_f} - \langle\varepsilon(\bu) \cdot \nu, &\bvarphi\rangle_{\Gamma_s} + (\nabla \pi, \bvarphi) \\
&= (\bu^{*}, \bvarphi) \quad \text{ for all }\,\bvarphi \in H_{\Gamma_f, 0}^{1}(\Omega_f),
\end{alignedat}
\label{maximality_calc_1}
\end{equation} 
and use of Green's theorem on the pressure term yields:
\begin{equation}
(\nabla \pi, \bvarphi)_{\Omega_f} = -(\pi, \div(\bvarphi))_{\Omega_f} + \langle \pi, \bvarphi \cdot \nu\rangle_{\Gamma_s}.
\label{maximality_calc_2}
\end{equation} 
Then combining \eqref{maximality_calc_1} with \eqref{maximality_calc_2} and the boundary coupling in (\hyperlink{domain4b}{D.4b})  $\varepsilon(\bu) \cdot \nu = \sigma(\bw) \cdot \nu + \pi\nu$ gives: 
\begin{equation}
\begin{alignedat}{2}
\lambda\,(\bu,\bvarphi)_{\Omega_f}
+(\varepsilon(\bu),\varepsilon(\bvarphi))_{\Omega_f}
- (\pi,\div\bvarphi)_{\Omega_f}&-\;\langle \sigma(\bw) \cdot \nu,\bvarphi\rangle_{\Gamma_s}
\\
&\;=\;(\bu^*,\bvarphi)_{\Omega_f} \quad \text{ for all }\,\bvarphi \in H_{\Gamma_f, 0}^{1}(\Omega_f),
\end{alignedat}
\label{maximality_calc_3}
\end{equation}
where for the first boundary term with $\sigma(\bw)\cdot\nu$, recall that \begin{equation}\begin{split}-\langle\sigma(
\bw) \cdot \nu, \bvarphi|_{\Gamma_s}\rangle_{\Gamma_s} &= (\sigma(\bw), \varepsilon(D_{\lambda}(\bvarphi|_{\Gamma_s})))_{\Omega_s} + (\div\sigma(\bw), D_{\lambda}(\bvarphi|_{\Gamma_s}))_{\Omega_s},
\\
&= (\sigma(\bw), \varepsilon(D_{\lambda}(\bvarphi|_{\Gamma_s})))_{\Omega_s} + ([\lambda^{2} + 1] \bw - (\lambda \bw^{*} + \bz^{*}), D_{\lambda}(\bvarphi|_{\Gamma_s}))_{\Omega_s}\end{split}\end{equation} where for the first equality we use the extension (right inverse of the trace) from $\bvarphi|_{\Gamma_s}$ to $\bvarphi \in \bH^{1}(\Omega_s)$, and for the second equality we recall the previous relation $\lambda^{2}\bw \; -\; \div\sigma(\bw) \;+\;\bw\; =\; \lambda \bw^{*} \; +\; \bz^{*}$ in $\Omega_s$ giving $\div\sigma(\bw) = (\lambda^{2} + 1)\bw - \lambda \bw^{*} - \bz^{*}$. We note that for first equality, the negative sign is using the fact that $\nu$ is interior to $\Omega_s$. Now plugging this back in to \eqref{maximality_calc_3} and also using $\bw = \bw_1 + \bw_0$ by \eqref{structure_via_fluid}, we get:
\begin{equation}
\begin{split}
\lambda\,(\bu,\bvarphi)_{\Omega_f}
&+(\varepsilon(\bu),\varepsilon(\bvarphi))_{\Omega_f}
- (\pi,\div\bvarphi)_{\Omega_f}
\\
&+\;(\sigma(\bw_1+\bw_0),\varepsilon(D_{\lambda}(\bvarphi|_{\Gamma_s})))_{\Omega_s}
\\
&+\;(([\lambda^{2} + 1](\bw_1+\bw_0)-\lambda\,\bw^{*}-\bz^{*},D_{\lambda}(\bvarphi|_{\Gamma_s}))_{\Omega_s}\;=\;(\bu^*,\bvarphi)_{\Omega_f},
\end{split}
\end{equation}
and then plugging in for $\bw_1 = \frac{1}{\lambda}D_{\lambda}(\bu|_{\Gamma_s})$ gives:
\begin{equation}
\begin{split}
\lambda\,(\bu,\bvarphi)_{\Omega_f}
&+(\varepsilon(\bu),\varepsilon(\bvarphi))_{\Omega_f}
- (\pi,\div\bvarphi)_{\Omega_f}
\\
&+\;(\sigma(\tfrac{1}{\lambda}D_{\lambda}(\bu|_{\Gamma_s})),\varepsilon(D_{\lambda}(\bvarphi|_{\Gamma_s})))_{\Omega_s}
+\;(\sigma(\bw_0),\varepsilon(D_{\lambda}(\bvarphi|_{\Gamma_s})))_{\Omega_s}
\\
&+\;([\lambda^{2} + 1]\,\tfrac{1}{\lambda}D_{\lambda}(\bu|_{\Gamma_s}),D_{\lambda}(\bvarphi|_{\Gamma_s}))_{\Omega_s}
+\;[\lambda^{2} + 1]\,(\bw_0,D_{\lambda}(\bvarphi|_{\Gamma_s}))_{\Omega_s}
\\
&-\;((\lambda\,\bw^*+\bz^*),D_{\lambda}(\bvarphi|_{\Gamma_s}))_{\Omega_s}
\;=\;(\bu^*,\bvarphi)_{\Omega_f}.
\end{split}
\end{equation}
Then we shift three data terms (containing $\bw_0$, $\bw^{*}$ and $\bz^{*}$) to the RHS to get:
\begin{equation}
\begin{split}
\lambda\,(\bu,\bvarphi)_{\Omega_f}
&+(\varepsilon(\bu),\varepsilon(\bvarphi))_{\Omega_f}
- (\pi,\div\bvarphi)_{\Omega_f}
\\
&+\;\tfrac{1}{\lambda}\,(\sigma(D_{\lambda}(\bu|_{\Gamma_s})),\varepsilon(D_{\lambda}(\bvarphi|_{\Gamma_s})))_{\Omega_s}
\\
&+\;\frac{\lambda^{2} + 1}{\lambda} \,(D_{\lambda}(\bu|_{\Gamma_s}),D_{\lambda}(\bvarphi|_{\Gamma_s}))_{\Omega_s}
=-\;(\sigma(\bw_0),\varepsilon(D_{\lambda}(\bvarphi|_{\Gamma_s})))_{\Omega_s}
\\
&
- [\lambda^{2} + 1]\,(\bw_0,D_{\lambda}(\bvarphi|_{\Gamma_s}))_{\Omega_s}
+ ((\lambda\,\bw^*+\bz^*),D_{\lambda}(\bvarphi|_{\Gamma_s}))_{\Omega_s}
+ (\bu^*,\bvarphi)_{\Omega_f}.
\end{split}
\end{equation}
Now we define $a( \cdot ; \cdot , \cdot)$ to be all of the terms on the LHS except for the $(\pi, \div\bvarphi)_{\Omega_f}$ which will be part of the bilinear form $b( \cdot, \cdot)$, and the forming term $F(\cdot)$ will be the terms in the RHS.

Specifically we define the bilinear form $a_{\lambda}(\cdot, \cdot): \bH_{\Gamma_f, 0}^{1}(\Omega_f) \times \bH_{\Gamma_f, 0}^{1}(\Omega_f) \to \mathbb{R}$ as:
\begin{equation}
\begin{alignedat}{2}
a_{\lambda}(\bv, \bvarphi)
  &\equiv \lambda(\bv, \bvarphi)_{\Omega_f} + (\varepsilon(\bv), \varepsilon(\bvarphi))_{\Omega_f} + \frac{1}{\lambda}(\sigma(D_{\lambda}(\bv|_{\Gamma_s})), \varepsilon(D_{\lambda}(\bvarphi|_{\Gamma_s})))_{\Omega_s} \\ &+ \frac{\lambda^{2} + 1}{\lambda}(D_{\lambda}(\bv|_{\Gamma_s}), D_{\lambda}(\bvarphi|_{\Gamma_s}))_{\Omega_s} 
  \qquad \qquad \text{ for all }\,\bv, \bvarphi \in \bH_{\Gamma_f, 0}^{1}(\Omega_f),
\end{alignedat}
\label{bilinear}
\end{equation}
with forcing term $F \in (\bH_{\Gamma_f, 0}^{1}(\Omega_f))^{*}$ defined by
\begin{equation}
\begin{alignedat}{2}
F(\bvarphi)
  &\equiv (\bu^{*}, \bvarphi)_{\Omega_f} + (\lambda \bw^{*} + \bz^{*}, D_{\lambda}(\bvarphi|_{\Gamma_s}))_{\Omega_s} \\& - \left(\sigma\left(\frac{1}{\lambda}D_{\lambda}(\bw^{*}|_{\Gamma_s}) + \mathcal{L}_{\lambda}^{-1}(\lambda \bw^{*} + \bz^{*})\right), \varepsilon(D_{\lambda}(\bvarphi|_{\Gamma_s}))\right)_{\Omega_s}  \\&- (\lambda^{2} + 1)\left(\frac{1}{\lambda}D_{\lambda}(\bw^{*}|_{\Gamma_s}) + \mathcal{L}_{\lambda}^{-1}(\lambda \bw^{*} + \bz^{*}), D_{\lambda}(\bvarphi|_{\Gamma_s})\right)_{\Omega_s}
  &\quad& \text{ for all }\,\bvarphi \in \bH_{\Gamma_f, 0}^{1}(\Omega_f).
\end{alignedat}
\label{mixed_forcing}
\end{equation}
In addition we define the bilinear form $b(\cdot, \cdot) : \bH_{\Gamma_f, 0}^{1}(\Omega_f) \times L^{2}(\Omega_f) \to \mathbb{R}$ by:
\begin{equation}
\begin{alignedat}{2}
b(\bvarphi, \mu)
  &\equiv -(\mu, \div\bvarphi)_{\Omega_f}
  &\quad& \text{ for all }\,\bvarphi \in \bH_{\Gamma_f, 0}^{1}(\Omega_f), \mu \in L^{2}(\Omega_f).
\end{alignedat}
\label{bilinear2}
\end{equation}
Then we have a mixed variational problem of finding a pair $[\bu, \pi] \in \bH_{\Gamma_f, 0}^{1}(\Omega_f) \times L^{2}(\Omega_f)$ that solves:
\begin{equation}
\begin{alignedat}{2}
a_{\lambda}(\bu,\bvarphi) \;+\; b(\bvarphi,\pi)
  &= F(\bvarphi)
  &\quad& \text{ for all }\,\bvarphi\in \bH_{\Gamma_f,0}^{1}(\Omega_f),\\
b(\bu,\mu)
  &= 0
  &\quad& \text{ for all }\mu\in L^{2}(\Omega_f).
\end{alignedat}
\label{mixed_variational}
\end{equation}

\noindent The subsequent steps follow that of \cite{avalos2008maximality}: verification of the inf-sup and other analytic properties, which establishes a unique solution $[\bu, \pi] \in \bH_{\Gamma_f, 0}^{1}(\Omega_f) \times L^{2}(\Omega_f)$ for the given $[\bu^{*}, \bw^{*}, \bz^{*}] \in \bH$. Then, by using the pair $[\bu, \pi]$ obtained to recover the solution components $\bw$ and $\bz$ and show that $[\bu, \bw, \bz] \in D(\mathcal{A})$ by the same procedure as in \cite{avalos2008maximality}, except with minor adjustments since we only use trace regularity for $\varepsilon(\bu) \cdot \nu - \pi\nu \in \bH^{-1/2}(\Gamma_s)$ (see Appendix~\ref{appendix:maximality_proof_2}). First we recall the Babu\v{s}ka--Brezzi theorem:

\begin{mytheorem}[Babu\v{s}ka--Brezzi, see e.g. p.116 of \cite{kesavan1989topics}]
\label{LBB}
Let $\bX$, $M$ be Hilbert spaces and $a : \bX \times \bX \to \mathbb{R}$, $b : \bX \times M \to \mathbb{R}$ continuous bilinear forms. Let
\begin{equation}\bZ := \{\bbeta \in \bX : b(\bbeta, \varrho) = 0 \quad \text{ for every } \varrho \in M \}.\end{equation}
\noindent Suppose that $a(\cdot, \cdot)$ is $\bZ$-elliptic, i.e. there exists a constant $\alpha > 0$ such that
\begin{equation}
a(\bbeta, \bbeta) \geq \alpha\|\bbeta\|_{\bX}^{2} \quad \text{ for every } \bbeta \in \bZ.
\end{equation}
Also suppose there exists a constant $\beta > 0$ such that
\begin{equation}
\sup_{\btau \in \bX}\frac{b(\btau, \varrho)}{\|\btau\|_{\bX} } \geq \beta \|\varrho\|_{M} \quad \text{ for every } \varrho \in M.
\label{infsup}
\end{equation}
Then for any $\boldsymbol{\kappa} \in \bX$ and $\ell \in M$, there exists a unique pair $(\hat{\bbeta}, \hat{\varrho}) \in \bX \times M$ such that:
\begin{equation}
\begin{alignedat}{2}
a(\hat{\bbeta},\btau) \;+\; b(\btau,\hat{\varrho})
  &= (\boldsymbol{\kappa}, \btau)_{\bX}
  &\quad& \text{ for all }\,\btau\in \bX,\\
b(\hat{\bbeta},\varrho)
  &= (\ell, \varrho)_{M}
  &\quad& \text{ for all }\varrho\in M.
\end{alignedat}
\label{mixed_variational_2}
\end{equation}
\end{mytheorem}  

Now we apply Theorem~\ref{LBB} for the system in \eqref{mixed_variational_0}, with $\bX = \bH_{\Gamma_f, 0}^{1}(\Omega_f)$, $M = L^{2}(\Omega_f)$, $(\boldsymbol{\kappa}, \bvarphi)_{\bX} = F(\bvarphi)$ for all $\bvarphi \in \bX$ and $(\ell, \mu) = 0$ for all $\mu \in M$, and verify conditions (i)--(iii) below: 

\begin{enumerate}[label=(\roman*)]

\item  $\bX$ and $M$ are Hilbert spaces and $a_{\lambda}( \cdot , \cdot)$ and $b( \cdot , \cdot )$ are bilinear continuous forms, which follow from the regularity of $D_{\lambda}$ and $\mathcal{L}_{\lambda}$ in \eqref{solid_BVP_1} and \eqref{solid_BVP_2}.

\item $a_{\lambda}(\cdot, \cdot)$ is elliptic over all of $\bX$, i.e. there exists a constant $\alpha > 0$ such that $a_{\lambda}(\bvarphi, \bvarphi) \geq \alpha \|\bvarphi\|_{\bX}^{2}$ for every $\bvarphi \in \bX$. Here by definition of $a_{\lambda}$ we have \begin{equation}a_{\lambda}(\bvarphi, \bvarphi) \geq  \|\varepsilon(\bvarphi)\|_{0, \Omega_f}^2 \qquad \text{ for all } \bvarphi \in \bX,\end{equation}
which gives ellipticity via Korn's inequality and Poincar\'e's inequality (using the Dirichlet boundary condition on $\Gamma_f$).

\item We have satisfaction of the `inf-sup' condition \eqref{infsup} for $b(\cdot, \cdot)$. This is exactly the conclusion established in Proposition~\ref{newproposition}, for $\bX = \bSigma$ and $M = L^{2}(\Omega_f)$.
\end{enumerate}

Hence by (i)--(ii), the Babu\v{s}ka--Brezzi Theorem can be applied to yield a unique pair $[\bu, \pi] \in \bX \times M = \bH_{\Gamma_f, 0}^{1}(\Omega_f) \times L^{2}(\Omega_f)$ which solves the system \eqref{mixed_variational}.

\subsection{Proof of maximality: Recovery of other variables}
\label{appendix:maximality_proof_2}

We get $[\bu, \pi] \in \bH_{\Gamma_f, 0}^{1}(\Omega_f) \times L^{2}(\Omega_f)$ solves \eqref{mixed_variational} and we infer from the second equation that:
\begin{equation}
\bu \in \bH_{\Gamma_f, 0}^{1}(\Omega_f)  \quad \text{ ; }  \quad \div(\bu) = 0 \text{ in } \Omega_f.
\end{equation}
We integrate by parts on the first equation to get:
\begin{equation}
\begin{alignedat}{2}
\lambda(\bu, \bvarphi)_{\Omega_f} - (\text{div}(\varepsilon(\bu)), \bvarphi)_{\Omega_f} + (\nabla \pi_0, \bvarphi)_{\Omega_f} &\quad& \text{ for all }\,\bvarphi \in \mathbf{\mathcal{D}}(\Omega_f).
\end{alignedat}
\end{equation}
Hence
\begin{equation}
\begin{alignedat}{2}
\lambda \bu - \text{div}(\varepsilon(\bu)) + \nabla \pi_0 = \bu^*,
\end{alignedat}
\end{equation}
where notice $\div(\bu) = 0$ and $\bu|_{\Gamma_f} = 0$ and hence (\hyperlink{domain1}{D.1}), so $\lambda \bu \in \mathcal{\bH}_f$ and so $-\div(\varepsilon(\bu)) + \nabla \pi_0 \in \mathcal{\bH}_f$ giving (\hyperlink{domain4a}{D.4a}). This gives $\varepsilon(\bu) \cdot \nu - \pi_{0}\nu \in \bH^{-1/2}(\Gamma_s)$.

\par Then we can recover the elastic variable $\bw$ by means of the relation in \eqref{structure_via_fluid}. Then by construction we have $\bw \in \bH^{1}(\Omega_s)$ and \eqref{structure_bvp}, the expression for which giving $\div(\varepsilon(\bw)) \in \bL^{2}(\Omega_s)$ and hence (\hyperlink{domain2}{D.2}).

\par Now we consider whether $\bw$ satisfies the boundary relation of (\hyperlink{domain4b}{D.4b}). By the second equation of \eqref{mixed_variational} and \eqref{structure_via_fluid}, we get
\begin{equation}
\begin{alignedat}{2}
\lambda(\bu, \bvarphi)_{\Omega_f} &+ (\varepsilon(\bu), \varepsilon(\bvarphi))_{\Omega_f} - (\pi, \text{div}\bvarphi)_{\Omega_f} \\ 
& + (\sigma(\bw), \varepsilon(D_{\lambda}(\bvarphi|_{\Gamma_s})))_{\Omega_s} + (\lambda^{2} + 1)(\bw, D_{\lambda}(\bvarphi|_{\Gamma_s}))_{\Omega_s} \\
&= (\bu^{*}, \bvarphi)_{\Omega_f} + (\lambda \bw^{*} + \bz^{*}, D_{\lambda}(\bvarphi|_{\Gamma_s}))_{\Omega_s} \quad& \text{ for all }\,\bvarphi \in \bH_{\Gamma_f, 0}^{1}(\Omega_f).
\end{alignedat}
\end{equation}
Then integrate by parts on $(\varepsilon(\bu), \varepsilon(\bvarphi))_{\Omega_f} - (\pi, \div(\bvarphi))_{\Omega_f} = (-\div(\varepsilon(\bu)) + \nabla \pi, \bvarphi)_{\Omega_f} + \langle \varepsilon(\bu) \cdot \nu - \pi \nu, \bvarphi\rangle_{\Gamma_s}$ using $\varepsilon(\bu) \cdot \nu - \pi\nu \in \bH^{-1/2}(\Gamma_s)$, and also the trace regularity for $\sigma(\bw) \cdot \nu \in \bH^{-1/2}(\Gamma_s)$, and we have
\begin{equation}
\begin{alignedat}{2}
\lambda(\bu, \bvarphi)_{\Omega_f} &+ (-\text{div}(\varepsilon(\bu)) + \nabla \pi, \bvarphi)_{\Omega_f} + \langle\varepsilon(\bu) \cdot \nu - \pi \nu, \bvarphi\rangle_{\Gamma_s} \\
& - (\text{div}\sigma(\bw), D_{\lambda}(\bvarphi|_{\Gamma_s}))_{\Omega_s} - \langle\sigma(\bw) \cdot \nu, \bvarphi|_{\Gamma_s}\rangle_{\Gamma_s} \\
& + (\lambda^{2} + 1)(\bw, D_{\lambda}(\bvarphi|_{\Gamma_s}))_{\Omega_s} \\
&= (\bu^{*}, \bvarphi)_{\Omega_f} + (\lambda \bw^{*} + \bz^{*}, D_{\lambda}(\bvarphi|_{\Gamma_s}))_{\Omega_s} \quad& \text{ for all }\,\bvarphi \in \bH_{\Gamma_f, 0}^{1}(\Omega_f).
\end{alignedat}
\end{equation}
Then by applying \eqref{eq:fluideq} and \eqref{structure_bvp} we get:
\begin{equation}
\begin{alignedat}{2}
\langle\varepsilon(\bu) \cdot \nu - \pi \nu, \bvarphi\rangle_{\Gamma_s} - \langle\sigma(\bw) \cdot \nu, \bvarphi|_{\Gamma_s}\rangle_{\Gamma_s} = 0 \quad& \text{ for all }\,\bvarphi \in \bH_{\Gamma_f, 0}^{1}(\Omega_f),
\end{alignedat}
\end{equation}
and we infer (using surjectivity of the Sobolev trace map) that:
\begin{equation}
[\varepsilon(\bu) \cdot \nu - \pi \nu]_{\Gamma_s} = \sigma(\bw) \cdot \nu \text{ in } \bH^{-1/2}(\Gamma_s).
\end{equation}
Finally set
\begin{equation}
\bz = \lambda \bw - \bw^{*} \in \bH^{1}(\Omega_s),
\end{equation}
and from the \eqref{structure_bvp} we get
\begin{equation}
\bz|_{\Gamma_s} = [\bu + \bw^{*}]_{\Gamma_s} - \bw^{*}|_{\Gamma_s} = \bu|_{\Gamma_s},
\end{equation}
and get (\hyperlink{domain3}{D.3}) and also get the two structure equations (using (\hyperlink{domain3}{D.3}) and \eqref{structure_bvp} ), hence get the total range condition \eqref{resolvent}, with solution $[\bu, \bw, \bz] \in D(\mathcal{A})$ as specified in (\hyperlink{domain1}{D.1})--(\hyperlink{domain4}{D.4}). This shows the maximality property for $\mathcal{A}$ (where notice that unlike \cite{avalos2008maximality} we only needed trace regularity for $\varepsilon(\bu) \cdot \nu - \pi\nu \in \bH^{-1/2}(\Gamma_s)$.)

\bibliographystyle{amsplain}
\bibliography{references}

\end{document}